\documentclass[draft]{amsart}
\usepackage{amsmath,amsfonts,amsrefs,latexsym}

\newtheorem{theorem}{Theorem}

\newtheorem{lemma}[theorem]{Lemma}

\theoremstyle{definition}

\begin{document}
\title{On  the Linearized Balescu-Lenard Equation}
\author[R. M. Strain]{Robert M. Strain}

\address{Department of Mathematics \\ Harvard University \\  
Cambridge \\ MA 02138 \\ USA} 
\email{strain at math.harvard.edu}

\keywords{Balescu-Lenard, collisional Kinetic Theory, Plasma Theory} 
\subjclass[2000]{Primary: 76P05; 
Secondary: 82B40, 82C40, 82D10}
\date{Completed: March 21, 2006, Revised: October 2, 2006}

\begin{abstract}
The Balescu-Lenard equation  from  plasma physics is widely considered to include a highly accurate correction to Landau's fundamental collision operator.   
Yet so far it has seen very little mathematical study.
We perform an extensive linearized analysis of this equation, 
which includes determining the asymptotic behavior of the new components of the linearized operator and establishing time decay rates for the linearized equation.   
\end{abstract}
\maketitle
\tableofcontents

\section{Introduction to the Balescu-Lenard Equation}

\thispagestyle{empty}

The Balescu-Lenard equation is a widely accepted kinetic equation which describes the dynamics of a spatially homogeneous plasma. It is 
\begin{equation}
\frac{\partial F}{\partial t} 
 = 
 \nabla \cdot \left\{\int_{\mathbb{R}^3} {\bf B}(v,v-v^*,\nabla F) \left\{F_*\nabla F -(\nabla F)_* F\right\} dv^* \right\}.
 \label{eq:blFULL}
\end{equation}
$F=F(t,v)$ is the velocity distribution, $F_*=F(t,v^*)$ and $(\nabla F)_*=\nabla_{v^*} F(t, v^*)$.  \eqref{eq:blFULL} is a correction to the spatially homogeneous Landau equation \eqref{landau}.  But in stark contrast to \eqref{landau},  
the kernel ${\bf B}=(B_{ij})$ introduces a strong nonlocal nonlinearity:
\begin{equation}
B_{ij}(v,v-v^*, \nabla F)
\equiv
\left(\frac{2n_0e^4}{m_e^2}\right)\int_{|k|\le k_0} \frac{k_ik_j}{|k|^4}\frac{ \delta(k\cdot[v-v^*])}{|\epsilon(k, k\cdot v, \nabla F)|^2}
 dk.
\label{eq:kernel}
\end{equation}
The main advantage of this operator is it's inclusion of the effects of Debye shielding, which we will see can lead to very different behavior.   Shielding is encoded into the collision kernel through the plasma dispersion function:
\begin{equation}
\epsilon(k, k\cdot v,\nabla F)
\equiv
1+\frac{1}{|k|^2}\left(\frac{v_e^2}{\lambda_e^2}\right)\lim_{\eta\downarrow 0}
\int_{\mathbb{R}^3}\frac{k\cdot \nabla F(u)}{k\cdot [v-u]-i\eta} du.  
\label{epsilondef}
\end{equation}
In \eqref{eq:blFULL}, \eqref{eq:kernel}, \eqref{epsilondef}, $t\ge 0$ and the velocity is $v=(v_1,v_2,v_3)\in\mathbb{R}^3$. 
 Also, the wavenumber $k=(k_1,k_2,k_3)\in\mathbb{R}^3$
is cut-off at
$0<k_0<\infty$.  The parameter, $k_0$,  represents the wavenumber beyond which collisions are no longer ``grazing collisions''. 
The constant 
$n_0$ is the average density, $e$ is the charge, $\lambda_e$ is the Debye length, $v_e$ is the thermal speed and $m_e$ is the mass.

The kinetic equation \eqref{eq:blFULL} was derived in a much more complicated form by
Bogoliubov \cite{MR0043001,MR0136381} in 1946.  Then, in 1960, Lenard \cite{MR0167274} showed how to write the equation explicitly in terms of the distribution, as \eqref{eq:blFULL}.  Lenard argued that this equation formally satisfies the expected physical properties including positivity of $F$, the standard conservation laws and the $H$-Theorem.  Further, Landau's equation \eqref{landau} was derived as an approximation to \eqref{eq:blFULL}.  Independently and in the same year, Balescu \cite{MR0128922} also derived  \eqref{eq:blFULL}, which is now commonly known as the Balescu-Lenard equation.  Due to the nonlocal nonlinearity in the collision kernel, there has always been an air of extreme difficulty surrounding this model.  As far as we know  there is only one other mathematically oriented paper on the subject \cite{MR0334783} from the early 1970's.   For more information on the physical background and relevance of the Balescu-Lenard equation see for instance \cite{MR0128922,MR1650315,hazeltineWaelbroeck,MR0167274,liboff1998,MR684990,montgomery,nicholson,MR0153257,MR1942465}.

The Landau equation, proposed by Landau in 1936, is one of the most fundamental partial differential equations in plasma physics; see several of the references just cited.  The spatially homogenous Landau equation takes the same form as \eqref{eq:blFULL}, save that the kernel ${\bf B}(v,v-v^*, \nabla F)$ is replaced by ${\bf b}(v-v^*)=(b_{ij})$ with
\begin{equation}
b_{ij}(v-v^*)
=
\frac{L}{|v-v^*|}
\left\{ \delta _{ij}-\frac{(v_i-v^*_i)(v_j-v^*_j)}{|v-v^*|^2}\right\} .  
\label{landau}
\end{equation}
Above $L$ is a parameter which is logarithmically divergent.  $L$ is proportional to
$$
L\sim\int \frac{d\chi}{\chi},
$$
which is divergent at both zero and infinity and therefore requires a cut-off resulting in the well known Coulomb logarithm.  The truncation at infinity is needed because Landau's equation does not model the effects of wide angle collisions, which is also the rational for $k_0$ in the Balescu-Lenard case.  However the cut-off near zero is needed because the Coulomb potential decreases very slowly at large distances.  

The inclusion of the effects of Debye shielding at large distances causes a rapid decrease.  
One of the key advantages of the Balescu-Lenard collision operator over the Landau collision operator is that it does not require a cut-off at small wavenumbers and can therefore model very precisely an electrically neutral plasma. On the other hand, we will argue that the inclusion of effects at small wave numbers makes  the difference between these two operators enormous.

Let's try to make this difference more precise.
Roughly speaking, one can see the Coulomb logarithm in the Balescu-Lenard kernel as follows
$$
L(\omega, v, \nabla F)=\left(\frac{2n_0e^4}{m_e^2}\right)\int_0^{k_0} \frac{d|k|}{|k| |\epsilon(|k|\omega,|k|\omega\cdot v, \nabla F)|^2}.
$$
Switching to spherical coordinates, $k=|k|\omega$, we can write \eqref{eq:kernel} as
$$
B_{ij}(v,v-v^*, \nabla F)
=
\int_{\omega\in S^2}L(\omega, v, \nabla F) \omega_i \omega_j \delta(\omega\cdot[v-v^*])
 d\omega.
$$
And if $L$ is a constant this is just the Landau kernel \eqref{landau}
because 
$$
\int_{\omega\in S^2} \omega_i \omega_j \delta(\omega\cdot [v-v^*])
 d\omega=\frac{\pi}{|v-v^*|}
\left\{ \delta _{ij}-\frac{(v_i-v^*_i)(v_j-v^*_j)}{|v-v^*|^2}\right\} .
$$
But one of the key results of our analysis shows that the Balescu-Lenard kernel \eqref{eq:kernel}  can be very far away from the Landau kernel \eqref{landau} in the following sense.

Consider the normalized steady state Maxwellian
$$
\mu(v)=(2\pi)^{-3/2}e^{-|v|^2/2}.
$$
We show  that up to some lower order decay, the kernel behaves like
$$
B_{ij}(v,v-v^*, \nabla \mu)  \approx \frac{ C}{|v-v^*|}
\frac{e^{\frac{1}{2}|v_R|^2}}{\left[1+|v_R|\right]^{3+\delta}}.
$$ 
Here $v_R$ is $v$ in the direction perpendicular to the relative velocity $v-v^*$:
\begin{equation}
|v_R|^2=|v|^2-\left(\frac{v-v^*}{|v-v^*|}\cdot v\right)^2.
\label{eq:vR}
\end{equation}
See Theorem \ref{lem:kernelestimate} for a precise statement.
The main new difficulty in our analysis is contained in this observation that the effect of Debye shielding on the Balescu-Lenard kernel,  when evaluated at maxwellian, is to create an exponentially growing velocity factor.  
Further, since \eqref{eq:blFULL} satisfies the H-Theorem, we speculate that this exponentially growing factor will be present for solutions to the Balescu-Lenard equation \eqref{eq:blFULL} at the nonlinear level for large times.
We consider $B_{ij}(v,v-v^*, \nabla \mu)$ above because this turns out to be the kernel of the linearized Balescu-Lenard collision operator.

Before stating our main results, we will linearize \eqref{eq:blFULL}.
For suitable functions $F$, $G$ and $H$ define the Balescu-Lenard collision operator by
\begin{equation*}
Q[F,G,H]
= 
 \nabla \cdot \left\{\int_{\mathbb{R}^3} {\bf B}(v,v-v^*,\nabla H) \left\{F_*\nabla G -(\nabla F)_*G \right\} dv^* \right\}.
\end{equation*}
We linearize this operator around the normalized steady state Maxwellian.
To this end, consider the standard perturbation 
$$
F=\mu+\sqrt{\mu}f.  
$$
Then we can write the Balescu-Lenard equation for the perturbation as
$$
\frac{\partial f}{\partial t} +Lf
 = 
N(f),
$$
where  the linearized Balescu-Lenard collision operator takes the form
\begin{equation}
Lf = -\mu^{-1/2}\left\{Q[\mu,\sqrt{\mu} f,\mu]+Q[\sqrt{\mu} f,\mu,\mu]  \right\}.
\label{linearizedBLop}
\end{equation}
The non-linear part of the Balescu-Lenard collision operator is
$$
N(f)
 =
\mu^{-1/2}\left\{Q[\mu+\sqrt{\mu}f, \mu+\sqrt{\mu}f, \mu+\sqrt{\mu}f]
-
Q[\mu,\sqrt{\mu} f,\mu]-Q[\sqrt{\mu} f,\mu,\mu]\right\}.
$$
At first glance, due to the way we have written it, $L$ may seem like a fabricated linear operator.  But  a Taylor expansion of the kernel ${\bf B}(v,v-v^*,\nabla\{\mu+\sqrt{\mu}f\})$ reveals that $N(f)$ is a nonlinear function of the perturbation $f$.  The terms subtracted off on the right cancel with the linear terms in the taylor expansion.  

More precisely, for $F$ satisfying $\epsilon(k, k\cdot v, \nabla F) \ne 0$, by \eqref{eq:kernel} we have 
$$
Q[\mu, \mu, F] = 0.
$$
Therefore,
$$
Q[\mu+\sqrt{\mu}f, \mu+\sqrt{\mu}f, F]
=
Q[\sqrt{\mu}f, \mu, F]
+
Q[\mu, \sqrt{\mu}f, F]
+
Q[\sqrt{\mu}f, \sqrt{\mu}f, F].
$$
Let us now look at the collision kernel \eqref{eq:kernel} with \eqref{epsilondef} for $F=\mu+\sqrt{\mu}f$.  From \eqref{epsilondef} 
$$
\epsilon(k, k\cdot v, \nabla \{\mu+\sqrt{\mu}f)
=
\epsilon(k, k\cdot v, \nabla \mu)
+
\frac{1}{|k|^2}\left(\frac{v_e^2}{\lambda_e^2}\right)\lim_{\eta\downarrow 0}
\int_{\mathbb{R}^3}\frac{k\cdot \nabla \{\sqrt{\mu}f\}(u)}{k\cdot [v-u]-i\eta} du.
$$
From here we see the first term in a taylor expansion of ${\bf B}(v,v-v^*,\nabla\{\mu+\sqrt{\mu}f\})$ in terms of the rightmost term above is  ${\bf B}(v,v-v^*,\nabla\mu)$ and all other terms depend on $f$.  Therefore 
\eqref{linearizedBLop} is the  linearized Balescu-Lenard collision operator.

Then the linearized Balescu-Lenard equation is
\begin{equation}
\frac{\partial f}{\partial t} +Lf=0.
\label{linearizedBL}
\end{equation}
The presence of physical constants does not create intrinsic mathematical difficulties.  Accordingly, to simplify our presentation, we will normalize all constants to one.
Because of the null space of $L$ (Lemma \ref{nullL}), a solution to \eqref{linearizedBL} formally satisfies
$$
\int_{\mathbb{R}^3} f_{0}(v)[1,v,|v|^{2}]\sqrt{\mu(v)} dv
=
\int_{\mathbb{R}^3} f(t,v)[1,v,|v|^{2}]\sqrt{\mu(v)} dv.
$$
We are interested in the asymptotic properties of solutions to \eqref{linearizedBL}.

We will prove time decay to maxwellian in weighted energy spaces.  
Consider the velocity weight
\begin{equation}
w(\ell ,\vartheta )(v)\equiv (1+|v|^{2})^{ \ell /2}
\exp \left( \frac{q}{4}(1+|v|^{2})^{\frac{\vartheta }{2}}\right).
\label{weight}
\end{equation}
Above $\ell\in\mathbb{R}$, $q>0$ and $0\le \vartheta\le 2$.  If $\vartheta=2$ we further assume $0<q<1$.  
 Depending on our choice of parameters, the velocity in this weight can grow either with an arbitrarily low polynomial power or alternatively almost but not quite as fast as $\mu^{-1/2}$.
Define the following weighted norm 
\begin{equation*}
|g|_{\vartheta}^2 \equiv \int_{\mathbb{R}^3} w^2(\ell ,\vartheta )|g(v)|^2 dv.
\end{equation*}
Above $\partial_i=\partial_{v_i}$. Let $\langle\cdot, \cdot \rangle$ denote the standard $L^2(\mathbb{R}^3_v)$ inner product. We put a zero in the norm to drop the entire weight \eqref{weight}, e.g. if $\ell=\vartheta=0$ then 
\begin{equation*}
|g|_{0}^2 = \int_{\mathbb{R}^3} |g|^2 dv=\langle g, g\rangle.
\end{equation*}
In these norms we can show linear decay:

\begin{theorem}\label{thm:MAIN}     Let $f_{0}(v)$ satisfy the conservation of mass, momentum and energy:
\begin{equation*}
\langle f_{0},[1,v,|v|^{2}]\sqrt{\mu }\rangle=0.
\end{equation*}
If $\vartheta>0$, a solution to \eqref{linearizedBL} with initial data $f(0,v)=f_0(v)$ and $|f_0|_{\vartheta}<\infty$ satisfies 
\begin{equation}
|f(t)|_{0}\le Ce^{-\lambda t^{p}}| f_0|_\vartheta,
\label{timeDECAY}
\end{equation}
where $C,\lambda>0$ and $p=\frac{\vartheta}{\vartheta+1}$.  
\end{theorem}

It is our hope that this linear decay will aid in a future in construction of classical solutions with small amplitude to the full Balescu-Lenard equation.

The only previous result for the the Balescu-Lenard equation that we know of is the work of Merchant and Liboff \cite{MR0334783} from 1973.  They show that the spectrum of the linearized Balescu-Lenard  operator is continuous from zero to minus infinity.  Additionally they obtain analytic expressions for some spherical harmonic eigenfunctions.  This lack of a spectral gap in the linear operator makes it difficult to prove time decay.  

Indeed, the proof of Theorem \ref{thm:MAIN} requires some development because there are many complicated elements of the Balescu-Lenard operator even at the linear level.   In Section \ref{sec:plasmaD} we will analyse the pointwise behavior of the longitudinal permittivity $\epsilon(k, k\cdot v, \nabla\mu)$.  We write down another formula for $\epsilon$ in \eqref{unperturbedplasmadispersion}.  Then we use this formula to determine the asymptotics of $\epsilon$ and thereby show that the Balescu-Lenard kernel $B_{ij}(v,v-v^*,\nabla\mu)$ is well defined (Lemma \ref{lem: Jdecay}).

Then, in Section \ref{sec:kernel}, we establish an alternate formula for the Balescu-Lenard kernel $B_{ij}(v,v-v^*, \nabla \mu)$ in Lemma \ref{lem:kernelREPRESENTATION} via a series of changes of variables.  We use this representation to show that the kernel contains an exponentially growing factor in Theorem \ref{lem:kernelestimate}.  In the proof of Theorem \ref{lem:kernelestimate} we split the integration region in order to squeeze a bit of extra decay out of the kernel, which we later need to use to estimate the linear operator.   

In Section \ref{sec:cf}, we look at the eigenvalues of the so-called ``collision frequency'', which is defined as
\begin{equation}
\label{sigma}
\sigma^{ij}(v)
\equiv
\int B_{ij}(v,v-v^*,\nabla \mu) \mu(v^*) dv^*.  
\end{equation}
Despite the exponential difference in the growth of the kernel as compared to the Landau kernel \eqref{landau}, the eigenvalues of the  collision frequency decay similarly to the Landau case.   Using a different series of changes of variables, we find a useful formula for the eigenvalues in and around Lemma \ref{lem:sigma}.  Then in Lemma \ref{lem:sigmaeigs}  we prove that the eigenvalues decay for large 
$|v|$ as
$$
\lambda_1(v) \sim   c_1 \log(2+|v|) [1+|v|]^{-3},
~~~~~~
\lambda_2(v) \sim   c_2[1+|v|]^{-1}.
$$
In the Landau case, the decay is as above if we remove the $\log(2+|v|)$ factor \cite{MR1463805}.

This study motivates the definition of the following weighted Sobolev norm
\begin{equation}
|g|_{\boldsymbol{\sigma},\vartheta}^2
\equiv
\sum_{i,j}\int_{\mathbb{R}^3}w^2(\ell ,\vartheta ) \sigma^{ij}(v) \left( \partial_i g \partial_j g + \frac{v_i}{2}\frac{v_j}{2}|g|^2\right) dv.  \label{normsigmaW}
\end{equation}
If $\vartheta$ is absent, this means that we drop the entire weight as follows
\begin{equation*}
|g|_{\boldsymbol{\sigma}}^2
=
|g|_{\boldsymbol{\sigma},0}^2
\equiv
\sum_{i,j}\int_{\mathbb{R}^3} \sigma^{ij}(v) \left( \partial_i g \partial_j g + \frac{v_i}{2}\frac{v_j}{2}|g|^2\right) dv.  
\end{equation*}
These turn out to be anisotropic spaces which are motivated by related spaces used by Guo in the Landau case \cite{MR1946444}.  They help measure the dissipation of \eqref{linearizedBLop}.

We study the full linear operator, $L$, in Section \ref{sec:linear}.  We write down it's null space in Lemma \ref{nullL}.  Then we split the linearized collision operator as
$$
\langle L f, f\rangle = |f |_\sigma^2 - \langle K f, f\rangle.
$$
See Lemma \ref{AKrepresent}.
The main estimate in our paper is the following

\begin{theorem}\label{thm:upper} For any small $\eta>0,$ there is $0<C(\eta)<\infty$ such that 
\begin{gather*}
|\langle w^2\partial_i \{\sigma^{ij}v_j\} g_1,g_2\rangle |  
+
|\langle w^2 Kg_1,g_2\rangle |  
\le 
\eta
|g_1|_{\boldsymbol{\sigma},\vartheta}|g_2|_{\boldsymbol{\sigma},\vartheta}
+
C(\eta)|g_1 {\bf 1}_{C(\eta)}|_{0}|g_2 {\bf 1}_{C(\eta)}|_{0}.
\end{gather*}
Here $w^2=w^2(\ell ,\vartheta)$ and ${\bf 1}_{C(\eta)}(v)$ is the indicator function for the set $\{|v|\le C(\eta)\}$.
\end{theorem}

Due to the exponential growth of the collision kernel (Theorem \ref{lem:kernelestimate}) it was difficult to expect that such a result was even true, especially with the possible inclusion of these exponentially growing weights.   However the exponential growth of the kernel is only present in the direction perpendicular to the relative velocity and there is some decay remaing in the other directions.   We design a splitting to show that this left over decay is just barely enough to prove Theorem \ref{thm:upper}.

In the rest of Section \ref{sec:linear}, we use Theorem \ref{thm:upper} to deduce coercivity of the linear operator \eqref{linearizedBLop} and exponential decay of solutions to the linearized equation \eqref{linearizedBL}.  
We use a standard compactness argument to prove the coercivity of the linear operator.  But all other arguments in this paper are constructive, including the estimate of $K$.  To prove Theorem \ref{thm:MAIN}, we use a splitting between velocity and time  which was previously used to establish decay rates for soft potential kinetic equations in \cite{MR575897,strainGUOed}.

To end this section, we remark that the results in each of the following sections build on and crucially make use of the results in previous sections.  We begin by analysing the longitudinal permittivity.

\section{The Dispersion function, $\epsilon(|k|,\hat{k}\cdot v,\nabla \mu)$}\label{sec:plasmaD}

The main objective of this section is to determine the pointwise asymptotic behavior of the plasma dispersion function \eqref{epsilondef} evaluated at maxwellian.  Since this appears to be difficult to accomplish from \eqref{epsilondef} due to the nonlocal operator our first step is to derive an alternate formula.  This formula 
\eqref{unperturbedplasmadispersion} seems to be known but we were unable to find a complete derivation in any one reference.  Therefore we  briefly derive \eqref{unperturbedplasmadispersion}  for completeness.    Then we establish the asymptotics 
 of $\epsilon$ in Lemma \ref{lem: Jdecay}.

Since we focus on the linear operator, which includes the dispersion function  only evaluated at Maxwellian
$\epsilon(k, k\cdot v,\nabla \mu)$,
in the rest of this article we will drop the dependence of the maxwellian  in our notation, to write only
$
\epsilon(k, k\cdot v).
$
Now we begin our computation of $\epsilon$ at Maxwellian.  From \eqref{epsilondef} we write
$$
\epsilon(k,s)
\equiv
\epsilon(k,s,\nabla \mu)
=
1+\frac{1}{|k|^2}\lim_{\eta\downarrow 0}
\int_{\mathbb{R}^3}\frac{k\cdot \nabla \mu(u)}{s-k\cdot u-i\eta} du.
$$
If $|k|\ne 0$, let $k^1=\frac{k}{|k|}$ and $\{k^1, k^2, k^3\}$ be an orthonormal basis for $\mathbb{R}^3$.  Define an orthogonal matrix $O_k$ such that $k\cdot O_ku=|k|u_1$ as 
\begin{equation}
\label{changeO}
O_k=[k^1~k^2~k^3].
\end{equation}
By this change of coordinates,
\begin{eqnarray*}
\int_{\mathbb{R}^3}\frac{k\cdot \nabla \mu(u)}{s-k\cdot u-i\eta} du
&=&
-(2\pi)^{-3/2}\int_{\mathbb{R}^3}\frac{(k\cdot u) e^{-|u|^2/2}}{s-k\cdot u-i\eta} du
\\
&=&
-(2\pi)^{-3/2}\int_{\mathbb{R}^3}\frac{ |k|u_1 e^{-|u|^2/2}}{s-|k|u_1-i\eta} du
\\
&=&
-(2\pi)^{-1/2}\int_{-\infty}^{+\infty}\frac{ |k|u_1 e^{-u_1^2/2}}{s-|k|u_1-i\eta} du_1.
\end{eqnarray*}
In the last step we integrated out the extra variables $u_2$ and $u_3$.  

To go further, we recall the well known Plemelj formula:
\begin{equation*}
\label{pjemel}
\lim_{\eta \downarrow 0} \frac{1}{x-y- i\eta} 
=
\text{P.V.}\left( \frac{1}{x-y} \right) + i\pi \delta(x-y),
\end{equation*}
which is a distriubtion.  Here ``P.V.'' denotes the Cauchy principle value.  Thus,
\begin{eqnarray*}
\lim_{\eta\downarrow 0}\int_{-\infty}^{+\infty}\frac{|k| u_1 e^{-u_1^2/2}}{s-|k|u_1-i\eta} du_1
&=&
\lim_{\eta\downarrow 0}\int_{-\infty}^{+\infty}\frac{ u_1e^{-u_1^2/2}}{\frac{s}{|k|}-u_1-i\frac{\eta}{|k|}} du_1
\\
&=&
\text{P.V.}\int_{-\infty}^{+\infty}\frac{ u_1 e^{-u_1^2/2}}{\frac{s}{|k|}-u_1} du_1
+i\frac{\pi s}{|k|}  e^{-\frac{s^2}{|k|^2} /2}.
\end{eqnarray*}
The second term on the r.h.s. is in exactly the form we want.  For the first term, 
further decompose $u_1=-\left(\frac{s}{|k|}-u_1\right)+\frac{s}{|k|}$ to obtain
\begin{eqnarray*}
\text{P.V.}\int_{-\infty}^{+\infty}\frac{ u_1 e^{-u_1^2/2}}{\frac{s}{|k|}-u_1} du_1
&=&
\frac{s}{|k|} \text{P.V.}\int_{-\infty}^{+\infty}\frac{  e^{-u_1^2/2}}{\frac{s}{|k|}-u_1} du_1
-
\int_{-\infty}^{+\infty}e^{-u_1^2/2} du_1,
\\
&=&
\frac{s}{|k|} \text{P.V.}\int_{-\infty}^{+\infty}\frac{  e^{-u_1^2/2}}{\frac{s}{|k|}-u_1} du_1
-
\sqrt{2\pi}.
\end{eqnarray*}
We collect the last few computations to obtain
$
\epsilon(|k|, x)
=
1+|k|^{-2}\Psi(x),
$
where 
 $\Psi=\Psi_R+i\Psi_I$ and
\begin{equation}
\begin{split}
\Psi_R(x)
&=
1-(2\pi)^{-1/2}x \text{P.V.}\int_{-\infty}^{+\infty}
\frac{  e^{-y^2/2}}{x-y} dy,
\\
\Psi_I(x)
&=
- \sqrt{\frac{\pi}{2}} x e^{-x^2/2}.
\end{split}
\label{psiI}
\end{equation}
The integral part of $\Psi_R$ is a well studied function in plasma physics \cite{MR0135916}.  

We further evaluate $\Psi_R$ via the well known formula \cite[p. 184-185, eq. (8.2.13)]{MR0153257}:
$$
\text{P.V.}\int_{-\infty}^{+\infty}\frac{ e^{- y^2/2}}{x-y} dy
=
\sqrt{2\pi} e^{-x^2/2}\int^{x}_{0}e^{y^2/2}dy.
$$
This yields the following simplified version of \eqref{psiI}
\begin{equation}
\label{psiR}
\Psi_R(x)
=
1-x e^{-x^2/2}\int^{x}_{0}e^{y^2/2}dy.
\end{equation}
By plugging $x=\frac{s}{|k|}=\frac{k}{|k|}\cdot v=\hat{k}\cdot v$ into $\Psi$, we can write
\begin{equation}
\label{unperturbedplasmadispersion}
\epsilon(|k|,\hat{k}\cdot v)
=
1+|k|^{-2}\Psi(\hat{k}\cdot v).
\end{equation}
Notice that above $\Psi$ only depends upon the magnitude of $v$ in the direction of $k$.  

Now that we have a tractable formula for $\epsilon$ evaluated at Maxwellian \eqref{unperturbedplasmadispersion}, we will study  it's pointwise behavior.  The key term in the decay of $\epsilon$ is $\Psi_R$:

\begin{lemma}\label{lem: Jdecay}  As $x\to\pm\infty$,
$
x^2\Psi_R(x)\to- 1.
$
\end{lemma}

\begin{proof}  We first write $x^2\Psi_R(x)$ as a fraction:
$$
x^2\Psi_R(x)
=x^2\left(1-xe^{-x^2/2}\int^{x}_{0}e^{y^2/2}dy\right)
=
\frac{\frac{1}{x} e^{x^2/2}-\int^{x}_{0}e^{y^2/2}dy}{\frac{1}{x^3} e^{x^2/2}}.
$$
Next, we apply l'H\^{o}pital's rule to obtain as $x\to\pm\infty$
$$
\frac{\frac{1}{x} e^{x^2/2}-\int^{x}_{0}e^{y^2/2}du}{\frac{1}{x^3} e^{x^2/2}}
\sim
\frac{ e^{x^2/2}\left(1-\frac{1}{x^2}\right)-e^{x^2/2}}{e^{x^2/2}\left(\frac{1}{x^2} -\frac{3}{x^4} \right)}
=
\frac{-1}{1 -\frac{3}{x^2} }\to -1.
$$
That's it.
\end{proof}

\noindent {\bf Remark \ref{lem: Jdecay}.1}.  
In particular, Lemma \ref{lem: Jdecay} shows that $\Psi_R(x)$, given by
\eqref{psiR},
is eventually negative for large enough $x$.  From (\ref{unperturbedplasmadispersion}) we have
$$
{\text{Re}}\{\epsilon(|k|,\hat{k}\cdot v)\}
=
1+\frac{\Psi_R(\hat{k}\cdot v)}{|k|^2}.
$$
Lemma \ref{lem: Jdecay} therefore tells us that the real part $\epsilon$ will be zero if $\hat{k}\cdot v$ is large enough 
and $|k|$ small enough--$|k|$ depending on the size of $\Psi_R(\hat{k}\cdot v)$.

There are then lots of places where $\text{Re}\{\epsilon\}$ is zero, but if $\text{Re}\{\epsilon(|k|,\hat{k}\cdot v)\}=0$ then fortunately
$$
{\text{Im}}\{\epsilon\}
=
\Psi_I
=
-\sqrt{\frac{\pi}{2}}\frac{(\hat{k}\cdot v)}{|k|^2}  e^{-(\hat{k}\cdot v)^2}\ne 0.
$$
We conclude that, when evaluated at maxwellian, $\epsilon\ne 0$ for all finite $v$.\\

By Lemma \ref{lem: Jdecay} and Remark \ref{lem: Jdecay}.1,  the kernel  of the Balescu-Lenard collision operator \eqref{eq:kernel} evaluated at Maxwellian is well defined.  In the next section we will further  use Lemma \ref{lem: Jdecay} to obtain asymptotic estimates of ${\bf B}(v,v-v^*,\nabla \mu)$ in Theorem \ref{lem:kernelestimate}.

\section{The Collision Kernel, ${\bf B}(v,v-v^*,\nabla \mu)$}\label{sec:kernel}

In this section we consider the Balescu-Lenard collision kernel \eqref{eq:kernel} at Maxwellian.  Since we will only consider this case, in the rest of the paper we will write
\begin{equation}
{\bf B}(v,v-v^*,\nabla \mu)={\bf B}(v,v-v^*).  
\label{kernATmax}
\end{equation}
We will look progressively more closely at the pointwise behavior of this kernel in a series of Lemmas.   The asymptotic analysis of the dispersion function from Section \ref{sec:plasmaD} will play an important role.

We first  record some basic properties of the collision kernel which are also shared with the Landau kernel \eqref{landau}.  A proof of the following can be found in \cite{montgomery}:
\begin{equation*}
{\bf B}\ge 0
~~ \text{and}~~
u^t {\bf B}(v,v-v^*) u=0
~~ 
\text{iff}
~~ 
u=c(v-v^*) 
~~
\text{with}
~~c\in \mathbb{R}.  
\end{equation*}
The main result of this section (in Theorem \ref{lem:kernelestimate}) is our pointwise estimate of the asymptotic growth rate  the collision kernel, which turns out to be exponential.  The first step in this direction is to develop a more tractable expression for \eqref{kernATmax}:

\begin{lemma}\label{lem:kernelREPRESENTATION}   
The Balescu-Lenard collision kernel \eqref{kernATmax} can be expressed as
\begin{equation*}
B_{ij}(v,v-v^*)
=
\frac{\xi_{ij}^1 w_1(|v_R|)+\xi_{ij}^2 w_2(|v_R|)}{|v-v^*|}.
\end{equation*}
Above, $\xi_{ij}^1$ and $\xi_{ij}^2$ are non-negative symmetric matricies that are connected to the Landau projection \eqref{landau} by the formula
$$
\xi_{ij}^1 +\xi_{ij}^2=\left(\delta_{ij}-\frac{(v_i-v^*_i) (v_j-v^*_j)}{|v-v^*|^2}\right).
$$
Further, 
\begin{equation}
\begin{split}
 \sum_{i} \xi_{ij}^1(v_i-v_i^*)= \sum_{i} \xi_{ij}^2(v_i-v_i^*)=0.
\end{split}
 \label{nullB}
\end{equation}
We define $w_1(|v_R|)$, $w_2(|v_R|)$ in \eqref{kernelWEIGHTS} and $\xi^1=(\xi_{ij}^1)$, $\xi^2=(\xi_{ij}^2)$  in \eqref{eq:xidef}.  
\end{lemma}

The weights are scalar functions which are given by
\begin{equation}
\begin{split}
w_1(|v_R|)
=
\int_0^{\pi/2}  \sin^2\theta
J(|v_R|\cos\theta)d\theta,
\\
w_2(|v_R|)
=
\int_0^{\pi/2}  \cos^2\theta
J(|v_R|\cos\theta)d\theta.
\label{kernelWEIGHTS}
\end{split}
\end{equation}
Here $J(s)$ includes the effects of the dispersion function \eqref{unperturbedplasmadispersion} as follows
\begin{equation}
J(x)
= 4\int_0^{k_0} 
 \frac{d\rho}{\rho |\epsilon(\rho, x)|^2}
= \int_0^{k_0} 
 \frac{4\rho^3}{\{\rho^2+\Psi_R( x )\}^2+\Psi_I^2( x )}d\rho.
 \label{eq:jay1}
\end{equation}
These quantities will help us get all our later estimates.

The main idea in the proof of Lemma \ref{lem:kernelREPRESENTATION} is to use three changes of coordinates, one at a time.  The first one is designed to extract the singularity from the delta function in \eqref{eq:kernel}.  The second coordinate change will give us a useful scalar quantity in the form of $|v_R|$.   And with the final rotation we obtain the integral $J$ which will be evaluated precisely in \eqref{lem:jayFORMULA} below.

\begin{proof}  Let $u^3=\widehat{v-v^*}=\frac{v-v^*}{|v-v^*|}$ and $\{u^1, u^2, u^3\}$ be an orthonormal basis.  Define $R$ to be the rotation matrix satisfying $Rk \cdot [v-v^*]=k_3 |v-v^*|$, e.g. 
$$
R=[u^1 \ u^2 \ u^3].
$$
Then we rotate the $k$ variable in \eqref{eq:kernel} 
with \eqref{unperturbedplasmadispersion}, in other words \eqref{kernATmax}, to obtain 
\begin{eqnarray*}
B_{ij}(v,v-v^*)
&\equiv&
\int_{|k|\le k_0} \frac{k_ik_j}{|k|^4}
\frac{ \delta(k\cdot[v-v^*])}{|\epsilon(|k|,\hat{k}\cdot v)|^2} dk
\\
&=&
\int_{|k|\le k_0} \frac{(Rk)_i(Rk)_j}{|k|^4}
\frac{ \delta(k_3|v-v^*|)}{|\epsilon(|k|,(R\hat{k})\cdot v)|^2} dk.
\end{eqnarray*}
With the formula $\delta(ak_3)=\frac{1}{a}\delta(k_3)$ for $a>0$ we have
\begin{eqnarray*}
B_{ij}(v,v-v^*)
&=&
\frac{1}{|v-v^*|}\int_{|k|\le k_0} \frac{(Rk)_i(Rk)_j}{|k|^4}
\frac{ \delta(k_3)}{|\epsilon(|k|,(R\hat{k})\cdot v)|^2} dk.
\end{eqnarray*}
We will expand $(Rk)_i(Rk)_j$ in order to simplify this expression.  
Note that
$$
(Rk)_i=k_1u_i^1+k_2u_i^2+k_3u_i^3.
$$
Here $u^l=(u^l_1, u^l_2, u^l_3)^t$ with $l\in\{1,2,3\}$.
We therefore have
$$
(Rk)_i(Rk)_j=\sum_{l=1}^3k_l^2u_i^l u_j^l+\sum_{l\ne m}k_l k_m u_i^l u_j^m.
$$
We  can therefore write
\begin{eqnarray*}
B_{ij}(v,v-v^*)
&=&
\frac{1}{|v-v^*|}
\int_{|k|\le k_0} \frac{\xi_{ij}}{|k|^2}
\frac{ \delta(k_3)}{|\epsilon(|k|,(R\hat{k})\cdot v)|^2} dk.
\end{eqnarray*}
Since terms involving $k_3$ vanish because of the delta function,
we can define
$$
\xi_{ij}
\equiv
|k|^{-2}\left\{k_1^2 u_i^1 u_j^1+k_2^2 u_i^2 u_j^2
+k_1 k_2\left(u_i^1 u_j^2 + u_j^1 u_i^2\right) \right\}.
$$
This completes the first change of variable.  With each new change of variable, below, we will elect to redefine $\xi_{ij}$ as needed instead of repeatedly introducing a new temporary notation.

Next we evaluate the delta function.  Consider $\underline{k}=(k_1,k_2)^t\in\mathbb{R}^2$ and write 
$$
\zeta=\frac{k_1}{|\underline{k}|} u^1+\frac{k_2}{|\underline{k}|} u^2.
$$
Evaluating the delta function yields
\begin{eqnarray*}
B_{ij}(v,v-v^*)
&=&
\frac{1}{|v-v^*|}
\int_{|\underline{k}|\le k_0} \frac{\xi_{ij}}{|\underline{k}|^2}
\frac{ 1}{|\epsilon(|\underline{k}|,\zeta\cdot v)|^2} d\underline{k},
\end{eqnarray*}
where we redefine
$$
\xi_{ij}
=
|\underline{k}|^{-2}\left\{k_1^2 u_i^1 u_j^1+k_2^2 u_i^2 u_j^2
+k_1 k_2\left(u_i^1 u_j^2 + u_j^1 u_i^2\right) \right\}.
$$
Before we rotate the $\underline{k}$  coordinate system again, we look at the following vector 
$$
v_R\equiv (v \cdot u^1,v \cdot u^2)^t.
$$
This is the velocity, $v$, in the direction perpendicular to the relative velocity, $v-v^*$.
Since $\{u^1, u^2, u^3\}$ is an orthonormal basis, we can expand
$$
|v|^2=(v\cdot u^1)^2+(v\cdot u^2)^2+(v\cdot u^3)^2.
$$
Thus,
$$
|v_R|^2=|v|^2-(v\cdot u^3)^2=|v|^2-\left(v\cdot \frac{v-v^*}{|v-v^*|}\right)^2.
$$
This is \eqref{eq:vR}.   Also
$$
\zeta \cdot v= \frac{\underline{k}}{|\underline{k}|}\cdot v_R.
$$
Now we are set up for another change of variables.

Let $O_R$ be the orthogonal matrix such that 
$
(O_R\underline{k})\cdot v_R=k_2 |v_R|,
$ 
e.g.
$$
O_R=[\hat{v}_R^\perp ~\hat{v}_R].
$$
Here 
$
v_R^\perp \equiv (u^2\cdot v, -u^1\cdot v)^t.
$  
We apply this rotation to obtain
\begin{eqnarray*}
B_{ij}(v,v-v^*)
&=&
\frac{1}{|v-v^*|}
\int_{|\underline{k}|\le k_0} \frac{\xi_{ij}}{|\underline{k}|^2}
\frac{ 1}{|\epsilon(|\underline{k}|, |v_R| k_2/ |\underline{k}|)|^2} d\underline{k},
\end{eqnarray*}
where $\xi_{ij}$ is redefined as
$$
\xi_{ij}
\equiv
|\underline{k}|^{-2}\left\{(O_R\underline{k})_1^2 u_i^1 u_j^1+(O_R\underline{k})_2^2 u_i^2 u_j^2
+(O_R\underline{k})_1 (O_R\underline{k})_2\left(u_i^1 u_j^2 + u_j^1 u_i^2\right) \right\}.
$$
To simplify this expression, we have to write out $O_R\underline{k}$.  
We expand
$$
\left( O_R \underline{k}\right)_i=k_1 \left( \hat{v}_R^\perp\right)_i+k_2 \left( \hat{v}_R \right)_i.
$$
Hence
\begin{eqnarray*}
\left( O_R \underline{k}\right)_i\left( O_R \underline{k}\right)_j
&=&
k_1^2 \left( \hat{v}_R^\perp \right)_i\left( \hat{v}_R^\perp \right)_j
+
k_2^2 \left( \hat{v}_R\right)_i  \left( \hat{v}_R\right)_j
\\
&&
+
k_1 k_2 \left\{\left( \hat{v}_R\right)_i \left( \hat{v}_R^\perp \right)_j+
\left( \hat{v}_R\right)_j \left( \hat{v}_R^\perp \right)_i \right\}.
\end{eqnarray*}
Plugging this into $\xi_{ij}$ yields a long  expression.  But the factor $k_1 k_2$ involves the integration of an odd function over an even domain, which is zero.  Disregarding terms with this factor, we can write 
$$
\xi_{ij}=\frac{k_1^2}{|\underline{k}|^2}\xi^1_{ij}+\frac{k_2^2}{|\underline{k}|^2}\xi^2_{ij},
$$
where
\begin{equation}
\begin{split}
\xi^1_{ij}
=&
\frac{(u^2\cdot v)^2 }{|v_R|^2} u^1_i u^1_j
+
\frac{(u^1\cdot v)^2 }{|v_R|^2} u^2_i u^2_j
-
\frac{(u^1\cdot v)(u^2\cdot v) }{|v_R|^2} 
\left(u_i^1 u_j^2 + u_j^1 u_i^2\right),
\\
\xi^2_{ij}
=&
\frac{(u^1\cdot v)^2 }{|v_R|^2} u^1_i u^1_j
+
\frac{(u^2\cdot v)^2 }{|v_R|^2} u^2_i u^2_j
+
\frac{(u^1\cdot v)(u^2\cdot v) }{|v_R|^2} \left(u_i^1 u_j^2 + u_j^1 u_i^2\right).
\end{split}
\label{eq:xidef}
\end{equation}
This completes our reduced expression after a second change of variables.

For the third and final change of variables, we will split into an angular integral and a magnitude integral.  The magnitude integral is evaluated precisely in \eqref{lem:jayFORMULA}.  For now, we choose polar coordinates as
$$
k_1=\rho \sin\theta,
~~
k_2=\rho \cos\theta.
$$
First changing coordinates and second plugging in \eqref{unperturbedplasmadispersion} yields
\begin{equation*}
\begin{split}
B_{ij}(v,v-v^*)
&=
\frac{1}{|v-v^*|}
\int_0^{2\pi}
\int_0^{k_0} 
\frac{ \left\{\xi_{ij}^1\sin^2\theta+\xi_{ij}^2\cos^2\theta  \right\} d\theta d\rho}
{\rho|\epsilon(\rho, |v_R| \cos\theta)|^2},
\\
&=
\frac{1}{|v-v^*|}
\int_0^{2\pi} 
\int_0^{k_0} 
\frac{ 1}{\rho} \frac{\left\{\xi_{ij}^1\sin^2\theta+\xi_{ij}^2\cos^2\theta  \right\} d\theta d\rho}{1+\frac{2\Psi_R}{\rho^2}+\rho^{-4}\left(\Psi_R^2+\Psi_I^2\right)},
\\
&=
\frac{1}{|v-v^*|}
\int_0^{2\pi} 
\left\{\xi_{ij}^1\sin^2\theta+\xi_{ij}^2\cos^2\theta  \right\} \frac{J(|v_R| \cos\theta)}{4}d\theta,
\end{split}
\end{equation*}
where the last line follows from the definition \eqref{eq:jay1}.

It remains to reduce the integral over $[0, 2\pi]$ to an integral over $[0,\pi/2]$.  Notice that $\Psi_R(x)$ is an even function \eqref{psiR} and  $\Psi_I(x)$ is odd \eqref{psiI}, so that   $\Psi_I^2(x)$ is even.  In particular, $J(x)=J(-x)$.  Since $J(x)$ is an even function \eqref{eq:jay1}, the reduction follows by first translating $\theta\to \pi-\theta$ on the region $[\pi, 2\pi]$ and second translating $\theta\to \frac{\pi}{2}-\theta$ on the region $[\pi/2, \pi]$.
\end{proof}

Next, it is not hard to evaluate \eqref{eq:jay1} precisely as:
\begin{equation}
J
=
 \log\left(1+\frac{k_0^4+2\Psi_R k_0^2}{\Psi_R^2+\Psi_I^2}\right)
+
2\frac{\Psi_R}{\Psi_I}
\left( \tan^{-1}\left(\frac{\Psi_R}{\Psi_I}\right)-\tan^{-1}\left(\frac{k_0^2+\Psi_R}{\Psi_I}\right)\right).
\label{lem:jayFORMULA}
\end{equation}
Here we clearly see the Logarithmic divergence of the Balescu-Lenard kernel when $k_0$ is sent to infinity.
With \eqref{lem:jayFORMULA}, we can determine the asymptotic limits of $J$:

\begin{lemma}\label{lem:jayLIMIT} 
$
x^3  e^{-x^2/2}J(x)\to \pm\sqrt{8\pi}
$
as
$
x\to\pm\infty.
$
\end{lemma}

Furthermore, if one were to add a cut-off at small wave number in the Balescu-Lenard kernel (say $k_s>0$) then it would appear in the first arctangent factor in \eqref{lem:jayFORMULA} as follows.  One would replace the factor $\frac{\Psi_R}{\Psi_I}$ by $\frac{k_s+\Psi_R}{\Psi_I}$. Our proof below implies this large exponential growth would then disappear because of cancellation.

\begin{proof} We will first examine the behavior of the log term.
Notice that the only zero of $ \Psi_I(x)$ is at $x=0$ and $ \Psi_R(0)=1$.  From this we see that
the argument of the log is bounded away from zero on any compact set.  So there are no finite singularities.  It is not hard to see that
as
$
x\to\pm\infty
$
$$
x^3  e^{-x^2/2}\log\left(1+\frac{k_0^4+2\Psi_R(x) k_0^2}{\Psi_R^2(x)+\Psi_I^2(x)}\right)\to 0.
$$
We therefore only need to look at the term involving the difference of tangents.

We now consider the second term of $J(x)$ in \eqref{lem:jayFORMULA}, which defines the asymptotics.
By Lemma \ref{lem: Jdecay} and \eqref{psiI}
$$
\frac{\Psi_R(x)}{\Psi_I(x)}
\sim
\frac{x^2\Psi_R(x)}{x^2\Psi_I(x)}
\sim
\frac{1}{x^3}\sqrt{\frac{2}{\pi}} e^{x^2/2}.
$$
Above ``$\sim$'' means the quantities have the same limit.  Similarly
$$
\frac{k_0^2+\Psi_R(x)}{\Psi_I(x)}
\sim
-\frac{1}{x}\sqrt{\frac{2}{\pi}} e^{x^2/2}.
$$
Since $\tan^{-1}(x)\to \pm \frac{1}{2}\pi$ as $x\to\pm\infty$, we conclude
$$
\lim_{x\to\pm\infty}\left( \tan^{-1}\left(\frac{\Psi_R(x)}{\Psi_I(x)}\right)-\tan^{-1}\left(\frac{k_0^2+\Psi_R(x)}{\Psi_I(x)}\right)\right)=\pm\pi.
$$
Therefore,
$$
\lim_{x\to\pm\infty}
x^3 e^{-x^2/2}2\frac{\Psi_R(s)}{\Psi_I(s)}
\left( \tan^{-1}\left(\frac{\Psi_R(s)}{\Psi_I(s)}\right)-\tan^{-1}\left(\frac{k_0^2+\Psi_R(s)}{\Psi_I(s)}\right)\right)
=\pm\sqrt{8\pi}.
$$
And the same is true for $x^3 e^{-x^2/2}J(x)$.  
\end{proof}

Given $u\in\mathbb{R}^3$, we define the projection
\begin{equation}
P_vu
\equiv 
(\hat{v}\cdot u)\hat{v}
=
\left(\frac{v}{|v|}\cdot u\right) \frac{v}{|v|}.
\label{vproj}
\end{equation}
Then in matrix form we have
\begin{gather*}
\begin{split}
(I-P_{v-v^*})_{ij}  &= \delta_{ij}-\frac{(v-v^*)_i (v-v^*)_j}{|v-v^*|^2}.
\end{split}
\end{gather*}
We now use Lemma \ref{lem:kernelREPRESENTATION} and Lemma \ref{lem:jayLIMIT}  to get asymptotic bounds for the kernel.

\begin{theorem}\label{lem:kernelestimate} Consider the Balescu-Lenard kernel \eqref{kernATmax} 
 and the relative velocity \eqref{eq:vR}.   For any $0<\delta<1$ there exists $C_\delta>0$ such that
\begin{equation*}
\left| B_{ij}(v,v-v^*) \right| 
\le
\frac{ e^{\frac{1}{2}|v_R|^2}}{|v-v^*|}
C_\delta\left[1+|v_R|\right]^{-3-\delta}.
\end{equation*}
Moreover, $\forall u\in\mathbb{R}^3$ and any $0<q<1$ there exists $C_q>0$ such that 
$$
u^t {\bf B}(v, v-v^*) u \ge \frac{|\{I-P_{v-v^*}\}u|^2}{C_q}\frac{e^{\frac{q}{2}|v_R|^2}}{|v-v^*|}
\left[1+|v_R|\right]^{-3}.
$$
In this sense, the upper bound is exponentially sharp.  
\end{theorem}

Although the bounds above are sufficient for the rest of our analysis, it seems tractable to refine the the bounds established in Theorem \ref{lem:kernelestimate} and make them sharp by evaluating more precisely the integrals in \eqref{intEVAL} below.

From Lemma \ref{lem:kernelREPRESENTATION}, \eqref{kernelWEIGHTS} and \eqref{eq:xidef} we see that
$$
\left| B_{ij}(v,v-v^*) \right| 
\le 4\frac{w_1(|v_R|)+w_2(|v_R|)}{|v-v^*|}.
$$
Therefore the upper bound in Theorem \ref{lem:kernelestimate} requires an upper bound on the weights $w_1$, $w_2$.  For the lower bound, we use \eqref{eq:xidef} to observe
\begin{equation*}
\begin{split}
u^t {\bf B}(v, v-v^*) u
&=\frac{w_1(|v_R|)}{|v-v^*|}u^t\xi^1u+\frac{w_2(|v_R|)}{|v-v^*|}u^t\xi^2u
\\
&\ge
\frac{\min\{w_1(|v_R|),w_2(|v_R|)\}}{|v-v^*|}\left(u^t\xi^1u+u^t\xi^2u\right).
\end{split}
\end{equation*}
Then using the formula for $\xi^1+\xi^2$ in Lemma \ref{lem:kernelREPRESENTATION} yields
$$
u^t\xi^1u+u^t\xi^2u=u^t(\xi^1+\xi^2)u=|\{I-P_{v-v^*}\}u|^2.
$$
We thereby see that for the lower bound in Theorem \ref{lem:kernelestimate} it is enough to get a lower bound on
$\min\{w_1(|v_R|),w_2(|v_R|)\}$.  This is what we prove.

\begin{proof}   We first establish the upper bound.  From \eqref{kernelWEIGHTS} and Lemma \ref{lem:jayLIMIT},  we have 
\begin{equation}
\begin{split}
w_1(|v_R|)
&\le 
C
\int_0^{\pi/2}  \sin^2\theta
\frac{e^{\frac{1}{2}|v_R|^2\cos^2\theta}}{[1+|v_R|\cos\theta]^3} d\theta,
\\
w_2(|v_R|)
&\le
C
\int_0^{\pi/2}  \cos^2\theta
\frac{e^{\frac{1}{2}|v_R|^2\cos^2\theta}}{[1+|v_R|\cos\theta]^3}
d\theta.
\label{intEVAL}
\end{split}
\end{equation}
Without loss of generality say $|v_R|\ge 1$.  
We split into $0\le \theta \le |v_R|^{-\delta}$ and $|v_R|^{-\delta}\le \theta \le \frac{\pi}{2}$.  
Then
\begin{eqnarray*}
\int_0^{|v_R|^{-\delta}}  
\frac{e^{\frac{1}{2}|v_R|^2\cos^2\theta}}{[1+|v_R|\cos\theta]^3}
d\theta
&\le&
\frac{1}{|v_R|^{\delta}}  
\frac{e^{\frac{1}{2}|v_R|^2}}{[1+|v_R|\cos(|v_R|^{-\delta} )]^3}.
\end{eqnarray*}
We will use 
$$
\cos(|v_R|^{-\delta} )\ge 1-\frac{1}{2!} |v_R|^{-2\delta}.
$$
From this lower bound and $0<\delta<1$ we have 
$$
1+|v_R|\cos(|v_R|^{-\delta}   )
\ge 1+|v_R|-\frac{1}{2}|v_R|^{1-2\delta}
\ge \frac{1}{2}(1+|v_R|).
$$
Here we have utilized $|v_R|\ge 1$.  Thus, 
$$
\frac{1}{|v_R|^{\delta}}  
\frac{e^{\frac{1}{2}|v_R|^2}}{[1+|v_R|\cos(|v_R|^{-\delta} )]^3}
\le 
2^{3+\delta}\frac{e^{\frac{1}{2}|v_R|^2}}{[1+|v_R|]^{3+\delta}}.
$$
This completes the estimate over  $0\le \theta \le |v_R|^{-\delta}$.

For the second half of the splitting, $|v_R|^{-\delta}\le \theta \le \frac{\pi}{2}$, we have
\begin{eqnarray*}
\int_{|v_R|^{-\delta}}^{\pi/2}
\frac{e^{\frac{1}{2}|v_R|^2\cos^2\theta}}{[1+|v_R|\cos\theta]^3}
d\theta
&\le&
C
e^{\frac{1}{2}|v_R|^2\cos^2(|v_R|^{-\delta})}.
\end{eqnarray*}
Now we will use the upper bound
$$
\cos(|v_R|^{-\delta})\le 1-\frac{1}{2!}|v_R|^{-2\delta}+\frac{1}{4!}|v_R|^{-4\delta}.
$$
From here, we get some weak exponential decay as long as $0<\delta<1$:
$$
e^{\frac{1}{2}|v_R|^2\cos^2(|v_R|^{-\delta})}
\le
C e^{\frac{1}{2}|v_R|^2}e^{-\frac{1}{4}|v_R|^{2(1-\delta)}}e^{\frac{1}{48}|v_R|^{2(1-2\delta)}}
\le
C e^{\frac{1}{2}|v_R|^2}e^{-\frac{1}{8}|v_R|^{2(1-\delta)}}.
$$  
This is more than enough decay to establish the upper bound.

Next we consider a lower bound for $\min\{w_1(|v_R|),w_2(|v_R|)\}$.    By Lemma \ref{lem:jayLIMIT},
\begin{equation*}
\begin{split}
w_1(|v_R|)
&\ge 
\frac{1}{C}
\int_0^{\pi/2}  \sin^2\theta
\frac{e^{\frac{1}{2}|v_R|^2\cos^2\theta}}{[1+|v_R|\cos\theta]^3} d\theta
\\
&\ge 
\frac{1}{C}[1+|v_R|]^{-3}
\int_0^{\pi/2}  \sin^2\theta
e^{\frac{1}{2}|v_R|^2\cos^2\theta}
d\theta,
\\
w_2(|v_R|)
&\ge
\frac{1}{C}
\int_0^{\pi/2}  \cos^2\theta
\frac{e^{\frac{1}{2}|v_R|^2\cos^2\theta}}{[1+|v_R|\cos\theta]^3}
d\theta
\\
&\ge
\frac{1}{C}[1+|v_R|]^{-3}
\int_0^{\pi/2}  \cos^2\theta
e^{\frac{1}{2}|v_R|^2\cos^2\theta}
d\theta.
\end{split}
\end{equation*}
This time we consider $w_1$ and $w_2$ separately.  For $w_2$ and $0<q<1$ we have
$$
\int_0^{\pi/2} \ge \int_0^{\cos^{-1}(q)} \cos^2\theta
e^{\frac{1}{2}|v_R|^2\cos^2\theta} d\theta
\ge q^2\cos^{-1}(q)  e^{\frac{q}{2}|v_R|^2}.
$$
And similarly for $w_1$,
$$
\int_0^{\pi/2} \ge \int_0^{\cos^{-1}(q)} \sin^2\theta
e^{\frac{1}{2}|v_R|^2\cos^2\theta} d\theta
\ge \left(\int_0^{\cos^{-1}(q)} \sin^2\theta d\theta\right) e^{\frac{q}{2}|v_R|^2}>0.
$$
These lower bounds for $w_1$ and $w_2$ establish the Theorem.
\end{proof}

\noindent {\bf Remark \ref{lem:kernelestimate}.1}.   It is a basic but important fact that ${\bf B}(v,v-v^*)={\bf B}(v^*,v-v^*)$, 
the delta function yields
\begin{eqnarray*}
B_{ij}(v,v-v^*)
&=&
\int_{|k|\le k_0} \frac{k_ik_j}{|k|^4}
\frac{ \delta(k\cdot[v-v^*])}{|\epsilon(|k|,\hat{k}\cdot v)|^2} dk
\\
&=&\int_{|k|\le k_0} \frac{k_ik_j}{|k|^4}
\frac{ \delta(k\cdot[v-v^*])}{|\epsilon(|k|,\hat{k}\cdot v^*)|^2} dk
=
B_{ij}(v^*,v-v^*).
\end{eqnarray*} 
We will use this later to split up the growth of $B_{ij}(v,v-v^*)$ in Theorem \ref{lem:kernelestimate} between $|v_R|$ and $|v_R^*|$.  Alternatively, the next identity can be used and additionally will be useful in other contexts below:
$$
|v_R|^2
=
{|v|^2-\left(\frac{(v-v^*)\cdot v}{|v-v^*|}\right)^2}
=
{|v^*|^2-\left(\frac{(v-v^*)\cdot v^*}{|v-v^*|}\right)^2}
=|v_R^*|^2.
$$
This is seen by a difference of squares argument:
$$
\left(\frac{(v-v^*)\cdot v}{|v-v^*|}\right)^2-\left(\frac{(v-v^*)\cdot v^*}{|v-v^*|}\right)^2
=
|v|^2-|v^*|^2.
$$
In particular, Theorem \ref{lem:kernelestimate} says that $B_{ij}(v, v-v^*)\in L^2_{loc}(\mathbb{R}^3_v\times \mathbb{R}^3_{v^*})$.

This completes our estimates for the collision kernel ${\bf B}(v,v-v^*)$.  In the next section we consider the  collision frequency.

\section{The Collision Frequency, $\sigma(v)$}\label{sec:cf}

We recall the Balescu-Lenard collision frequency \eqref{sigma}.
 We will use the reduced form of ${\bf B}$ from
Lemma \ref{lem:kernelREPRESENTATION} to do an asymptotic analysis of the collision frequency.  The following quantity will also be used:
\begin{equation}
\label{sigmai}
\sigma^{i}(v)
 \equiv 
\left\{\sigma(v)\frac{v}{2}\right\}_i
=\sum_{j} \int B_{ij}(v,v-v^*) \frac{v_j^*}{2} \mu(v^*) dv^*.
\end{equation}
In this section, we compute the eigenvalues $\lambda_1(v)$ and $\lambda_2(v)$ of $\sigma(v)$ in Lemma \ref{lem:sigma}.  Then in Lemma \ref{lem:sigmaeigs} we transform these eigenvalues into a form which is conducive to obtaining precise pointwise information.  Somewhat surprisingly, despite the exponential growth of the collision kernel (Theorem \ref{lem:kernelestimate}) in comparision to the Landau operator \eqref{landau}, the collision frequency still decays at large velocities, and in a way which is different but also closely related to the Landau case.

We also used the collision frequency to define an anisotropic norm 
\eqref{normsigmaW}, 
$|\cdot|_{\boldsymbol{\sigma}}$, which measures the dissipation of the linear operator $L$.  An equivalent norm will be established in Corollary \ref{lem:sigmaeigs}.1, using all the asymptotics developed in this section.  This norm is important for establishing the main results in Section \ref{sec:linear}.  First we look at the eigenvalues:

\begin{lemma}\label{lem:sigma}  
The matrix $\sigma(v)$ from \eqref{sigma} has an eigenvalue $\lambda_1(v)$ with eigenvector $v$ and a double eigenvalue $\lambda_2(v)$ whose eigenspace is perpendicular to $v$:
\begin{equation*}
\begin{split}
\lambda_1(|v|) 
& = 
(2\pi)^{-1/2}\int_{|k|\le k_0} \frac{k_1^2}{|k|^5}
\frac{ e^{-\frac{1}{2}(|v|k_1/|k|)^2}}{|\epsilon(|k|,k_1 |v|/|k|)|^2} dk,
\\
\lambda_2(|v|) 
& =  
(2\pi)^{-1/2}\int_{|k|\le k_0} \frac{k_2^2+k_3^2}{2|k|^5}
\frac{ e^{-\frac{1}{2}(|v|k_1/|k|)^2}}{|\epsilon(|k|,k_1 |v|/|k|)|^2} dk.
\end{split}
\end{equation*}
And we can expand
$
\sigma^{ij}(v)
=
\hat{v}_i \hat{v}_j \lambda_1(|v|)
+
\left(\delta_{ij}-\hat{v}_i \hat{v}_j\right) \lambda_2(|v|)
$
with 
$
\hat{v}=\frac{v}{|v|}.
$
\end{lemma}

We remark that the strategy for computing these eigenvalues is to use a series of two rotations in a different way from Lemma \ref{lem:kernelREPRESENTATION}.  Here we first evaluate the $v^*$ integration which is not present in $B_{ij}(v,v-v^*)$.

\begin{proof} First we translate $v^*\to v-v^*$ and second we use \eqref{changeO},  $O_k$,  to obtain
\begin{eqnarray*}
\sigma^{ij}(v) & = & \int B_{ij}(v,v-v^*)\mu(v^*)  dv^*
 = 
\iint_{|k|\le k_0} \frac{k_ik_j}{|k|^4}
\frac{ \delta(k\cdot[v-v^*])\mu(v^*)}{|\epsilon(|k|,\hat{k}\cdot v)|^2} dk dv^*,
 \\
 & = & 
\iint_{|k|\le k_0} \frac{k_ik_j}{|k|^4}
\frac{ \delta(k\cdot v^*)\mu(v-v^*)}{|\epsilon(|k|,\hat{k}\cdot v)|^2} dv^* dk,
 \\
 & = & 
\iint_{|k|\le k_0} \frac{k_ik_j}{|k|^4}
\frac{ \delta(k\cdot O_k v^*)\mu(v-O_k v^*)}{|\epsilon(|k|,\hat{k}\cdot v)|^2} dv^* dk,
 \\
 & = & 
\iint_{|k|\le k_0} \frac{k_ik_j}{|k|^5}
\frac{ \delta(v^*_1)\mu(v-O_k v^*)}{|\epsilon(|k|,\hat{k}\cdot v)|^2} dv^* dk.
\end{eqnarray*}
Next, with $k^1=\hat{k}$ and $\{k^1, k^2, k^3\}$ an orthonormal basis for $\mathbb{R}^3$ as in \eqref{changeO}, we expand the exponent of the Maxwellian as
$$
|v-O_k v^*|^2=|v-v_1^*k^1-v_2^*k^2-v_3^*k^3|^2.
$$
We can further write $v=(v\cdot k^1)k^1+(v\cdot k^2)k^2+(v\cdot k^3)k^3$.  Then by orthogonality
$$
\mu(v-Ov^*)=(2\pi)^{-3/2} e^{-\frac{1}{2}\left\{[(v\cdot k^1)-v_1^*]^2+[(v\cdot k^2)-v_2^*]^2+[(v\cdot k^3)-v_3^*]^2 \right\}}.
$$
Now we evaluate the delta function and use translation invariance to obtain
\begin{equation}
\begin{split}
\label{sigmak}
\sigma^{ij}(v)
&=
\iint_{|k|\le k_0} \frac{k_ik_j}{|k|^5}
\frac{ \delta(v^*_1)\mu(v-O_k v^*)}{|\epsilon(|k|,\hat{k}\cdot v)|^2} dv^* dk,
\\
 &= 
(2\pi)^{-1/2}  \int_{|k|\le k_0} \frac{k_ik_j}{|k|^5}
\frac{ e^{-\frac{1}{2}(v\cdot \hat{k})^2}}{|\epsilon(|k|,\hat{k}\cdot v)|^2} dk.
\end{split}
\end{equation} 
We will use this formula to compute the eigenvalues.

Now that we have evaluated the $v^*$ integrations via a rotation in the $k$ direction, we will simplify the $k$ integrations by rotating in the $v$ direction.   Let $v^1=\hat{v}$ and $\{v^1, v^2, v^3\}$ be an orthonormal basis for $\mathbb{R}^3$.  Further define the rotation matrix
\begin{equation}
\label{changeU}
O_v=[ v^1~ v^2~ v^3 ].
\end{equation}
Notice that
$
(O_v k)_i=k_1 v^1_i+k_2 v^2_i+k_3 v^3_i.
$
We rotate the $k$ variable with $O_v$ in \eqref{sigmak} to achieve
\begin{eqnarray*}
\sigma^{ij}(v) 
 & = & 
(2\pi)^{-1/2} \int_{|k|\le k_0} \frac{(O_vk)_i(O_vk)_j}{|k|^5}
\frac{ e^{-\frac{1}{2}(v\cdot O_v\hat{k})^2}}{|\epsilon(|k|,O_v\hat{k}\cdot v)|^2} dk,
\\
 & = & 
(2\pi)^{-1/2}
\int_{|k|\le k_0} \frac{ k_1^2v^1_i v^1_j+k_2^2v^2_i v^2_j+k_3^2 v^3_i v^3_j}{|k|^5}
\frac{e^{-\frac{1}{2}\left(\frac{k_1}{|k|}|v|\right)^2}}{|\epsilon(|k|, k_1|v|/|k|)|^2} dk.
\end{eqnarray*}
Above, the cross terms in $(O_vk)_i(O_vk)_j$ disappear because they give you the integral of an odd function over an even domain.  
By symmetry,
$$
\int_{|k|\le k_0} \frac{ 
k_2^2}{|k|^5}\frac{e^{-\frac{1}{2}\left(\frac{k_1}{|k|}|v|\right)^2}}{|\epsilon^m(|k|, k_1|v|/|k|)|^2} dk
=
\int_{|k|\le k_0} \frac{ 
k_3^2}{|k|^5}\frac{e^{-\frac{1}{2}\left(\frac{k_1}{|k|}|v|\right)^2}}{|\epsilon^m(|k|, k_1|v|/|k|)|^2} dk.
$$
Recall \eqref{vproj}.  By the spectral theorem
$
I-P_{v^1}=P_{v^2}+P_{v^3},
$
or in another form 
$$
v^2_i v^2_j +v^3_i v^3_j=\delta_{ij}-\hat{v}_i \hat{v}_j.
$$  
These last three points yield the result.  \end{proof}

With a sequence of two rotations, we can write these eigenvalues in the form
\begin{eqnarray*}
\lambda_1(v)
&=&
\iint_{|k|\le k_0}\left(\hat{k}\cdot \hat{v}\right)^2
\frac{ \delta(k\cdot v^*)\mu(v-v^*)}{|k|^2|\epsilon(|k|,\hat{k}\cdot v)|^2} dk dv^*,
\\
\lambda_2(v)
&=&
\iint_{|k|\le k_0}\frac{1}{2}\left|\frac{k}{|k|}\times\frac{v}{|v|}\right|^2
\frac{ \delta(k\cdot v^*)\mu(v-v^*)}{|k|^2|\epsilon(|k|,\hat{k}\cdot v)|^2} dk dv^*.
\end{eqnarray*}
This follows by first rotating the $v^*$ variable with \eqref{changeO} to obtain
\begin{eqnarray*}
\lambda_1(v)
&=&
(2\pi)^{-1/2} \int_{|k|\le k_0}
\frac{ \left(\hat{k}\cdot \hat{v}\right)^2e^{-\frac{1}{2}(\hat{k}\cdot v)^2}}{|k|^3|\epsilon(|k|,\hat{k}\cdot v)|^2} dk,
\\
\lambda_2(v)
&=&
(2\pi)^{-1/2}\int_{|k|\le k_0}\frac{1}{2}\left|\frac{k}{|k|}\times\frac{v}{|v|}\right|^2
\frac{ e^{-\frac{1}{2}(\hat{k}\cdot v)^2}}{|k|^3|\epsilon(|k|,\hat{k}\cdot v)|^2} dk.
\end{eqnarray*}
Next rotate the $k$ variable in the direction of $v$ with \eqref{changeU} to get  the eigenvalues as written in Lemma \ref{lem:sigma}.    These forms of the eigenvalues are similar in form to the eigenvalues found for the collision frequency of the Landau equation in \cite{MR1463805}.  

We will now  transform these eigenvalues into a form which is conducive to asymptotic analysis in Lemma \ref{lem:sigmaeigs} below.  We switch to spherical coordinates  via 
\begin{gather*}
k_1  =  \rho \cos\theta_1, ~~
k_2  =  \rho \sin\theta_1 \cos\theta_2, ~~
k_3  =  \rho \sin\theta_1 \sin\theta_2, 
\\
0\le \rho\le k_0, ~~ 0\le \theta_1\le \pi, ~~ 0\le \theta_2\le 2\pi.
\end{gather*}
Then we can write the eigenvalues in Lemma \ref{lem:sigma}  as
\begin{eqnarray*}
\lambda_1(|v|) 
& = &
\sqrt{2\pi}
\int_{0}^\pi d\theta_1  \int_0^{k_0} d\rho
 \frac{\cos^2\theta_1 \sin\theta_1}{\rho}
\frac{ e^{-\frac{1}{2}|v|^2\cos^2\theta_1}}{|\epsilon(\rho,|v|\cos\theta_1)|^2},
\\
\lambda_2(|v|) 
& = & 
\sqrt{\frac{\pi}{2}}
\int_{0}^\pi d\theta_1 \int_0^{k_0} d\rho
 \frac{\sin^3\theta_1}{\rho}
\frac{ e^{-\frac{1}{2}|v|^2\cos^2\theta_1}}{|\epsilon(\rho,|v|\cos\theta_1)|^2}.
\end{eqnarray*}
Plug \eqref{eq:jay1} into the eigenvalues to obtain
\begin{eqnarray*}
\lambda_1(|v|) 
& = &
\sqrt{\frac{\pi}{8}}
\int_{0}^\pi d\theta_1  
\cos^2\theta_1 \sin\theta_1
 e^{-\frac{1}{2}|v|^2\cos^2\theta_1}J(|v|\cos\theta_1),
\\
\lambda_2(|v|) 
& = & 
\frac{1}{4}\sqrt{\frac{\pi}{2}}
\int_{0}^\pi d\theta_1 
\sin^3\theta_1e^{-\frac{1}{2}|v|^2\cos^2\theta_1}J(|v|\cos\theta_1).
\end{eqnarray*}
We change variables as $y=|v|\cos\theta_1$ above to obtain
\begin{eqnarray*}
\lambda_1(|v|)
&=&
\frac{1}{|v|^3}\sqrt{\frac{\pi}{8}} \int_{-|v|}^{+|v|}
 y^2 e^{-\frac{1}{2}y^2}  J(y) dy,
\\
\lambda_2(|v|)
&=&
\frac{1}{|v|}\frac{1}{4}\sqrt{\frac{\pi}{2}} \int_{-|v|}^{+|v|}  \left(1-\frac{y^2}{|v|^2}\right) e^{-\frac{1}{2}y^2}  J(y) dy.
\end{eqnarray*}
The fact that the integrand's are even functions yields
\begin{eqnarray*}
\lambda_1(|v|)
&=&
 \frac{1}{|v|^3}\sqrt{\frac{\pi}{2}} \int_{0}^{|v|}
 y^2  e^{-\frac{1}{2}y^2}J(y) dy,
\\
\lambda_2(|v|)
&=&
\frac{1}{|v|}\sqrt{\frac{\pi}{8}} \int_{0}^{|v|}  \left(1-\frac{y^2}{|v|^2}\right) e^{-\frac{1}{2}y^2}J(y) dy.
\end{eqnarray*}
Now we are ready to look at the decay of these eigenvalues:

\begin{lemma}\label{lem:sigmaeigs}  
As $|v|\to\infty$
\begin{equation*}
\begin{split}
\lambda_1(|v|)& \sim  2\pi \frac{\log(2+|v|)}{1+|v|^3},
\\
\lambda_2(|v|)& \sim   \sqrt{\frac{\pi}{8}}\frac{ \int_{0}^{\infty}  e^{-\frac{1}{2}y^2}J(y) dy}{1+|v|}.
\end{split}
\end{equation*}
Moreover, for any multiindex $\beta$ with $|\beta|\le 1$ we have
\begin{equation}
\begin{split}
\label{destlambda}
\left| D^\beta \lambda_1(|v|)\right| & \le   C^1 \log(2+|v|)(1+|v|)^{-3-|\beta|}, 
\\
\left| D^\beta \lambda_2(|v|)\right|& \le  C^2 (2+|v|)^{-1-|\beta|},
\end{split}
\end{equation}
where $C^1$ and $C^2$ are positive constants.
\end{lemma} 

\noindent {\bf Remark \ref{lem:sigmaeigs}.1}.  As $|v|\downarrow 0$ these eigenvalues converge to a unique finite limit:
\begin{equation*}
\lim_{x\downarrow 0} \lambda_1(x)=\lim_{x\downarrow 0} \lambda_2(x)=\frac{1}{3}\sqrt{\frac{\pi}{2}}  J(0).
\end{equation*}
$J(0)$ is clearly finite by \eqref{eq:jay1}.  
The eigenvalues thus have no finite singularites.

\begin{proof}  
By Lemma \ref{lem:jayLIMIT} and l'H\^{o}pital's rule, as $|v|\to+\infty$
$$
 \frac{ \int_{0}^{|v|}
 y^2e^{-\frac{1}{2}y^2}  J(y) dy}{\log(|v|)}
 \sim
 |v|^3e^{-\frac{1}{2}|v|^2}  J(|v|) \to \sqrt{8\pi},
$$
which implies the decay for $\lambda_1$ in Lemma \ref{lem:sigmaeigs}.  

For the decay of $\lambda_2$, we consider 
$$
|v| \lambda_2(|v|)
=
\sqrt{\frac{\pi}{8}} \int_{0}^{|v|}  e^{-\frac{1}{2}y^2}J(y) dy
-\frac{1 }{|v|}\sqrt{\frac{\pi}{8}}\int_{0}^{|v|} y^2  e^{-\frac{1}{2}y^2}J(y) dy.
$$
Since the second term on the r.h.s. converges to zero, as $|v|\to+\infty$ we have
$$
|v| \lambda_2(|v|)
\to
\sqrt{\frac{\pi}{8}} \int_{0}^{\infty}  e^{-\frac{1}{2}y^2}J(y) dy.
$$
And the integral on the r.h.s. is finite by Lemma \ref{lem:jayLIMIT}.  This yields the decay of $\lambda_2$.  

The decay of $\lambda_1$ in \eqref{destlambda} follows from taking the derivative of $\lambda_1$ and then applying Lemma \ref{lem:jayLIMIT} and l'H\^{o}pital's rule, exactly as we have done for the no derivative case.  The same recipe will establish the decay of the derivative of $\lambda_2$ in \eqref{destlambda}.
\end{proof}

Then all the results in this section yield the following lower bound:

\medskip \noindent {\bf Corollary \ref{lem:sigmaeigs}.2}.
{\it $\exists c>0$ such that
\begin{gather*}
\frac{|g|_{\boldsymbol{\sigma},\vartheta}^2}{c}
\ge 
\int_{\mathbb{R}^3}w^2\left\{
\frac{\log(2+|v|)}{1+|v|^3}|P_v\nabla g|^2
+\frac{|[I-P_v] \nabla g|^2}{1+|v|}
+\frac{\log(2+|v|)}{1+|v|}|g|^2 \right\}dv,
\end{gather*}
where $w^2=w^2(\ell ,\vartheta )(v)$ is defined in \eqref{weight} and $P_v$ is defined in \eqref{vproj}.}

\begin{proof}  By \eqref{vproj}, Lemma \ref{lem:sigma}   and \eqref{normsigmaW},
$$
|g|_{\boldsymbol{\sigma},\vartheta}^2
=
\int_{\mathbb{R}^3}w^2(\ell ,\vartheta )\left\{\lambda_1(v) \left( |P_v\nabla g|^2 +|v|^2|g|^2\right)
+\lambda_2(v)  |[I-P_v]\nabla g|^2\right\} dv.
$$
Plugging in the asymptotics from Lemma \ref{lem:sigmaeigs} yields Corollary \ref{lem:sigmaeigs}.2.
\end{proof}

An analogous upper bound can also be established in the same way. 
This completes our study of the Balescu-Lenard collision frequency \eqref{sigma}.  In Section \ref{sec:linear} we will use these asymptotics to prove bounds for the Linearized Balescu-Lenard Operator.

\section{Compactness of $K$, Coercivity of $L$ and Exponential Decay}\label{sec:linear}

The main result of this section (Theorem \ref{thm:upper}) is to show that $L$ can be split into as $L=-A-K$ (Lemma \ref{AKrepresent}) where $K$ is ``compact''  (in the sense of the inequality in Theorem \ref{thm:upper}).
This is a standard Theorem for a linearized collision operator of a kinetic equation, such as the Boltzmann or the Landau equation.  However the exponential growth of the kernel 
(Theorem \ref{lem:kernelestimate}), which is not present in the Boltzmann or Landau kernel, creates new difficulties.  As a consequence of Theorem \ref{thm:upper}, we deduce coercivity for $L$ in the anisotropic norm \eqref{normsigmaW} in Corollary \ref{thm:upper}.1.  Then we finish  section \ref{sec:linear} by proving exponential decay of solutions to the linearized Balescu-Lenard equation.

First define the projection 
$$
{\bf P}g=\{a_g+{\bf b}_g\cdot v+c_g|v|^2\}\mu^{1/2}(v),
$$
where $a_g,c_g\in \mathbb{R}$ and ${\bf b}_g\in \mathbb{R}^3$ depend on the function $g(v)$.
Then we have the following standard lemma for the linearized collision operator \eqref{linearizedBLop}:

\begin{lemma}\label{nullL}  $\langle Lg_1, g_2\rangle = \langle g_1, Lg_2\rangle$, $L\ge 0$ and $\langle Lg_1,g_1\rangle=0$ if and only if $g_1(v)={\bf P}g_1$.  In this case $L{\bf P}g_1= 0$.  
\end{lemma}

For the Landau equation, there is a standard argument used to prove this Lemma, see for instance  \cite[Lemma 4]{MR1946444}.   The proof in the Balescu-Lenard case is exactly the same as for the Landau equation because the null space \eqref{nullB} is exactly the same.

Now we will split the operator $L=-A -K$ using \eqref{sigma} and \eqref{sigmai}:

\begin{lemma}\label{AKrepresent}  We split $L=-A -K$, where
\begin{eqnarray*}
Ag&=&
\partial_i \left\{\sigma^{ij} \partial_jg \right\}
+
\{ \partial_i \sigma^{i}\} g
-
\sigma^{ij}\frac{v_i}{2} \frac{v_j}{2}g,
\\
Kg
&=&
-\mu^{-1/2}\partial_i\left\{ \mu \int_{\mathbb{R}^3} B_{ij}(v,v-v^*)\mu^{1/2}_* 
\left\{(\partial_j g)_*+ \frac{v_j^*}{2}g_*\right\} dv^*\right\}.
\end{eqnarray*}
Here and in the proof we use the convention of summing over repeated indicies.
\end{lemma}

The proof of Lemma \ref{AKrepresent} is virtually the same as \cite[Lemma 1]{MR1946444}.  This is expected because the Landau kernel and the Balescu-Lenard kernel share the same null space.  

\begin{proof}     We define
\begin{gather*}
Ag = \mu^{-1/2}Q[\mu,\mu^{1/2} g,\mu],
\\
Kg=\mu^{-1/2}Q[\mu^{1/2} g,\mu,\mu]. 
\end{gather*}
Then $L=-A-K$ by \eqref{linearizedBLop}.  We will simplify $A$ and $K$.  
Notice that
\begin{gather*}
\partial_i \mu=-v_i \mu,
\\
\partial_i[\mu^{1/2} g]=\mu^{1/2}\left(\partial_i-\frac{v_i}{2}\right)g, 
\\
\partial_i[\mu^{-1/2} g]=\mu^{-1/2}\left(\partial_i+\frac{v_i}{2}\right)g.
\end{gather*}
We will use these and the null space  \eqref{nullB} of  \eqref{eq:kernel}  several times below. 

We compute
\begin{eqnarray*}
Ag
& = &  
\mu^{-1/2}
\partial_i \int_{\mathbb{R}^3} B_{ij}(v,v-v^*) 
\left\{\mu_* \partial_j(\mu^{1/2}g) -(\mu^{1/2}g) (\partial_j \mu)_* \right\} dv^*,
\\
 & = & 
\mu^{-1/2}
\partial_i \int_{\mathbb{R}^3} B_{ij}(v,v-v^*)\mu_*\mu^{1/2} 
\left\{\left(\partial_j-\frac{v_j}{2}\right)g+v^*_jg  \right\} dv^*,
\\
 & = & 
\mu^{-1/2}
\partial_i \int_{\mathbb{R}^3} B_{ij}(v,v-v^*)\mu_*\mu^{1/2} 
\left(\partial_j+\frac{v_j^*}{2}\right)g  dv^*,
\end{eqnarray*}
where we used \eqref{nullB} in the last step.  By \eqref{nullB} again, and then by \eqref{sigma} we have
\begin{eqnarray*}
 & = & 
\mu^{-1/2}
\partial_i \int_{\mathbb{R}^3} B_{ij}(v,v-v^*)\mu_*\mu^{1/2} 
\left(\partial_j+\frac{v_j}{2}\right)g  dv^*,
 \\
  & = & 
\mu^{-1/2}
\partial_i \left\{\sigma^{ij}\mu^{1/2} \left(\partial_j g+\frac{v_j}{2} g\right) \right\}.
\end{eqnarray*}
Next, we take the derivatives on each term to obtain
\begin{eqnarray*}
  & = & 
\mu^{-1/2}
\partial_i \left\{\mu^{1/2}\sigma^{ij} \partial_j g \right\}
+
\mu^{-1/2}
\partial_i \left\{\mu^{1/2}\sigma^{ij} \frac{v_j}{2} g \right\},
 \\
   & = & 
\partial_i \left\{\sigma^{ij} \partial_j g \right\}
-
\sigma^{ij}\frac{v_i}{2} \partial_j g 
+
\left\{\sigma^{ij} \frac{v_j}{2}  \right\}\partial_i g
+
\partial_i\left\{\sigma^{ij} \frac{v_j}{2}  \right\} g
-\sigma^{ij} \frac{v_i}{2}\frac{v_j}{2} g,
 \\
   & = & 
\partial_i \left\{\sigma^{ij} \partial_j g \right\}
+
\partial_i\left\{\sigma^{ij} \frac{v_j}{2}  \right\} g
-
\sigma^{ij} \frac{v_i}{2}\frac{v_j}{2} g.
\end{eqnarray*}
This is the desired expression for $A$.  

Next, for $K$, we have
\begin{eqnarray*}
Kg
& = &  
\mu^{-1/2}\partial_i \left\{\int_{\mathbb{R}^3} B_{ij}(v,v-v^*) 
\left\{(\mu^{1/2}g)_*\partial_j\mu -(\partial_j [\mu^{1/2}g])_* \mu \right\} dv^* \right\},
\\
 & = & 
\mu^{-1/2}\partial_i\left\{ \mu \int_{\mathbb{R}^3} B_{ij}(v,v-v^*)\mu^{1/2}_* \left\{-v_jg_* +\frac{v_j^*}{2}g_*-(\partial_j g)_* \right\} dv^*\right\}.
\end{eqnarray*}
By \eqref{nullB} this is 
\begin{eqnarray*}
 & = & 
-\mu^{-1/2}\partial_i\left\{ \mu \int_{\mathbb{R}^3} B_{ij}(v,v-v^*)\mu^{1/2}_* 
\left\{(\partial_j g)_*+ \frac{v_j}{2}g_*\right\} dv^*\right\},
\\
 & = & 
-\mu^{-1/2}\partial_i\left\{ \mu \int_{\mathbb{R}^3} B_{ij}(v,v-v^*)\mu^{1/2}_* 
\left\{(\partial_j g)_*+ \frac{v_j^*}{2}g_*\right\} dv^*\right\}.
\end{eqnarray*}
And this is the expression we sought for $K$.  
\end{proof}

We are ready to prove Theorem \ref{thm:upper}. \\

\noindent {\it Proof of Theorem \ref{thm:upper}.}
We first estimate $|\langle w^2(\ell ,\vartheta ) Kg_1,g_2\rangle|$.  Recall Lemma \ref{AKrepresent} and write
\begin{gather*}
\langle w^2(\ell ,\vartheta )Kg_1, g_2\rangle
\\
=
\int_{\mathbb{R}^3\times \mathbb{R}^3} w^2(\ell ,\vartheta )B_{ij}(v,v-v^*)\sqrt{\mu_*\mu} 
\left\{(\partial_jg_1)_*+\left(\frac{v_j}{2}g_1\right)_* \right\}
\left\{\partial_ig_2+\frac{v_i}{2}g_2 \right\} dv^* dv,
\\
+\int_{\mathbb{R}^3\times \mathbb{R}^3} \partial_i w^2(\ell ,\vartheta )B_{ij}(v,v-v^*)\sqrt{\mu_*\mu} 
\left\{(\partial_jg_1)_*+\left(\frac{v_j}{2}g_1\right)_* \right\}
g_2 dv^* dv,
\end{gather*}
where we recall that $(\cdot)_*$ means that the function in parenthesis is evaluated at $v^*$.  We will split the integration region several times to obtain the estimate.    The derivative of the weight function 
\eqref{weight} is 
\begin{gather}
\partial _i(w^2(\ell,\vartheta))=w^2(\ell,\vartheta)w_1v_i,
\label{wDERIVATIVE}
\end{gather}
where 
$$
w_1(v)=\left\{2\ell(1+|v|^2)^{-1}+q\frac{\vartheta}{2}(1+|v|^2)^{\frac{%
\vartheta}{2}-1}\right\}.  
$$
In particular $|w_1(v)|\le C<\infty$. Then we can write
\begin{gather}
\langle w^2Kg_1, g_2\rangle
=
\sum_{i,j}\int_{\mathbb{R}^3\times \mathbb{R}^3} w^2(\ell ,\vartheta )B_{ij}(v,v-v^*)\sqrt{\mu_*\mu} 
(d_1^jg_1)_*
d_2^i g_2 dv^* dv,
\label{K12}
\end{gather}
where
\begin{gather}
d_1^j
=
\partial_j+\frac{v_j}{2}, 
~~~
d_2^i
=
\partial_i+\frac{v_i}{2}+w_1(v)v_i, 
\label{hopdef}
\end{gather}
We will estimate $\langle w^2(\ell ,\vartheta )Kg_1, g_2\rangle$ in this condensed form.

Now we outline the main strategy of the proof.
A key point is to get sufficient upper bounds for $w(\ell,\vartheta)(v)B_{ij}(v,v-v^*)\sqrt{\mu_*\mu}$. 
We first want to control $B_{ij}(v,v-v^*)\sqrt{\mu_*\mu}$ by something which approximates the dissipation $\sigma^{ij}$ in our norm.  But we also want velocity decay left over to generate a small constant factor as in the statement of Theorem \ref{thm:upper}.  And further we want to show that this bound has additional velocity decay   to allow us to distribute and control the exponential weight \eqref{weight}, which only depends on $v$.  This is done via  several splittings.

We will now look for a bound for $w(\ell,\vartheta)(v)$ in terms of $w(\ell,\vartheta)(v^*)$.
To this end, we expand
$$
|v|^2=|v|^2-\left(\frac{v-v^*}{|v-v^*|}\cdot v\right)^2+\left(\frac{v-v^*}{|v-v^*|}\cdot v\right)^2.
$$
Then, since $0\le \vartheta\le 2$, we obtain
$$
e^{ \frac{q}{4}(1+|v|^{2})^{\frac{\vartheta }{2}}}
\le \exp \left( \frac{q}{4}\left(1+|v|^2-\left(\frac{v-v^*}{|v-v^*|}\cdot v\right)^2\right)^{\frac{\vartheta }{2}}\right)
e^{ \frac{q}{4}\left(\frac{v-v^*}{|v-v^*|}\cdot v\right)^{\vartheta}}.
$$
If $\ell\ge 0$ then we similarly have
$$
\left(1+|v|^{2}\right)^{ \ell /2}
\le  
C\left(1+|v|^2-\left(\frac{v-v^*}{|v-v^*|}\cdot v\right)^2\right)^{ \ell /2}\left(1+\left(\frac{v-v^*}{|v-v^*|}\cdot v\right)^2\right)^{ \ell /2}.
$$
But if $\ell<0$ then
$$
\left(1+|v|^{2}\right)^{ \ell /2}
\le  
\left(1+|v|^2-\left(\frac{v-v^*}{|v-v^*|}\cdot v\right)^2\right)^{ \ell /2}.
$$
Then, by the last few inequalities, under any conditions we have shown 
$$
w(\ell ,\vartheta )(v) \le Cw(\ell ,\vartheta )(v_R)
\left(1+\left(\frac{v-v^*}{|v-v^*|}\cdot v\right)^2\right)^{ |\ell| /2}
\exp \left( \frac{q}{4}\left(\frac{v-v^*}{|v-v^*|}\cdot v\right)^{\vartheta}\right).
$$
Since $|v_R|=|v_R^*|$ (Remark \ref{lem:kernelestimate}.1) we have
$$
w(\ell ,\vartheta )(v_R)=w(\ell ,\vartheta )(v_R^*)
\le 
Cw(\ell ,\vartheta )(v^*)\left(1+\left(\frac{v-v^*}{|v-v^*|}\cdot v\right)^2\right)^{ |\ell| /2}.
$$
The extra factor on the end is needed only if $\ell<0$.  We have this shown that
\begin{gather}
\label{weightDISTRIBUTE}
w(\ell ,\vartheta )(v)\le  w(\ell ,\vartheta )(v^*)
\left(1+\left(\frac{v-v^*}{|v-v^*|}\cdot v\right)^2\right)^{ |\ell|}
\exp \left( \frac{q}{4}\left(\frac{v-v^*}{|v-v^*|}\cdot v\right)^{\vartheta}\right).
\end{gather}
This estimate allows us to distribute the  exponentially growing velocity weight from the $v$ variable onto the $v^*$ variable.  We have to pay with some extra growth of $v$ in the direction of the relative velocity, but this can be controlled by terms in the upper bound \eqref{bijbound0} below.

The next step is to get bounds for $\left|w(\ell ,\vartheta )(v)B_{ij}(v, v-v^*) \sqrt{\mu_*\mu} \right|$.
  From  Theorem \ref{lem:kernelestimate} and Remark \ref{lem:kernelestimate}.1,
\begin{gather}
\label{bijbound0}
\left|B_{ij}(v, v-v^*) \sqrt{\mu_*\mu} \right|
\le
\frac{C\exp\left(-\frac{1}{4}\left(\frac{v-v^*}{|v-v^*|}\cdot v\right)^2
-\frac{1}{4}\left(\frac{v-v^*}{|v-v^*|}\cdot v^*\right)^2\right)}{|v-v^*|[1+|v_R|]^{3+\delta}}.
\end{gather}
So we have some exponential decay in the direction of the relative velocity, and this is how we control the exponentially growing weight \eqref{weight} and \eqref{weightDISTRIBUTE}.
Since either $0\le \vartheta<2$ and $0<q$ or $\vartheta=2$ and  $0<q<1$ there is $0<q^\prime<1$ such that
\begin{gather}
\nonumber
w(\ell ,\vartheta )(v)\left|B_{ij}(v, v-v^*) \sqrt{\mu_*\mu} \right|
\\
\le
\frac{Cw(\ell ,\vartheta )(v^*)}{|v-v^*|[1+|v_R|]^{3+\delta}}
\exp\left(-\frac{q^\prime}{4}\left(\frac{v-v^*}{|v-v^*|}\cdot v\right)^2
-\frac{1}{4}\left(\frac{v-v^*}{|v-v^*|}\cdot v^*\right)^2\right).
\label{bijbound00}
\end{gather}
We will use this upper bound several times below.

We first split the integration into a compact region and a large region.  For large $m_1>0$, define a smooth function $\varphi_1(r)$ so that $\varphi_1(r)=0$ for $r\le m_1$ and $\varphi_1(r)=1$ for 
$r\ge 2m_1$.  Then, with \eqref{K12} and \eqref{hopdef}, we define $K_1$  so that 
\begin{gather*}
\langle w^2K_1g_1, g_2\rangle
=
\sum_{i,j}\int_{\mathbb{R}^3\times \mathbb{R}^3} \varphi_1(|v|+|v^*|)w^2B_{ij}(v,v-v^*)\sqrt{\mu_*\mu} 
(d_1^jg_1)_*
d_2^i g_2 dv^* dv.
\end{gather*}
This is really the hardest term to estimate.  
We will split the integration region a few more times to do it.

We first split into two regions where \eqref{bijbound00} yields solid exponential decay in $|v|$ and $|v^*|$.  Define 
 $$
S_1=\{ |v|\ge 2 |v^*|   \}\cup \{  |v^*|\ge 2 |v|  \}.
$$ 
On $S_1$, if $|v|\ge 2 |v^*|$ then
$$
|v-v^*|\le |v|+|v^*|\le \frac{3}{2}|v|.
$$
And further
$$
\frac{(v-v^*)\cdot v}{|v-v^*|}
 \ge
\frac{|v|^2-|v^*||v|}{|v-v^*|}
 \ge
\frac{|v|^2-\frac{1}{2}|v|^2}{|v-v^*|}
=
\frac{1}{2}\frac{|v|^2}{|v-v^*|}
  \ge 
 \frac{1}{3} |v|
 \ge 
 \frac{2}{3} |v^*|.
$$
Alternatively, if $|v^*|\ge 2 |v|$ then
$$
|v-v^*|\le |v|+|v^*|\le \frac{3}{2}|v^*|.
$$
And similarly
\begin{gather*}
 \frac{(v^*-v)\cdot v^* }{|v-v^*|}
  \ge
 \frac{|v^*|^2-|v||v^*| }{|v-v^*|}
 \ge
 \frac{1}{2}\frac{|v^*|^2 }{|v-v^*|}
   \ge 
 \frac{1}{3} |v^*|
 \ge 
 \frac{2}{3} |v|.
\end{gather*}
In either case, we plug these estimates into \eqref{bijbound00} to obtain
\begin{gather}
\left|w(\ell ,\vartheta )(v)B_{ij}(v, v-v^*) \mu^{1/2}(v)\mu^{1/2}(v^*) \right|
\le
w(\ell ,\vartheta )(v^*)\frac{Ce^{-\frac{q^\prime}{36}(|v|^2+|v^*|^2)}}{|v-v^*|}.
\label{bijbound1}
\end{gather}
This is the strongest decay estimate we'll get.  

Here and below we will  define terms like $\langle {\bf 1}_{S_1} w^2K_{1}g_1, g_2\rangle$ to be
$$
\langle {\bf 1}_{S_1} w^2K_{1}g_1, g_2\rangle
=
\sum_{i,j}\int_{S_1} \varphi_1(|v|+|v^*|)w^2B_{ij}(v,v-v^*)\sqrt{\mu_*\mu} 
(d_1^jg_1)_*
d_2^i g_2 dv^* dv.
$$
On $S_1$ we use \eqref{bijbound1} (and recall \eqref{hopdef}) to get
\begin{gather*}
\left|\langle {\bf 1}_{S_1} w^2K_{1}g_1, g_2\rangle\right|
\le
C\int_{|v|+|v^*|\ge m_1} \frac{e^{-\frac{q^\prime}{36}(|v|^{2}+|v^*|^{2})}}{|v-v^*|}
\left|(w(\ell ,\vartheta )d_1^jg_1)_* \right|
\left|w(\ell ,\vartheta )d_2^i g_2  \right|
\\
\le
Ce^{-\frac{q^\prime}{144}m_1^{2}}
\int \frac{e^{-\frac{q^\prime}{72}(|v|^{2}+|v^*|^{2})}}{|v-v^*|}
\left|(w(\ell ,\vartheta )d_1^jg_1)_* \right|
\left|w(\ell ,\vartheta )d_2^i g_2  \right|
dv^* dv.
\end{gather*}
By Cauchy-Schwartz, this is
\begin{gather*}
\le
Ce^{-\frac{q^\prime}{144}m_1^{2}}
\left(\int \frac{e^{-\frac{q^\prime}{72}(|v|^{2}+|v^*|^{2})}}{|v-v^*|}
\left|(w(\ell ,\vartheta )d_1^jg_1)_* \right|^2dv^* dv\right)^{1/2}
\\
\times
\left(\int \frac{e^{-\frac{q^\prime}{72}(|v|^{2}+|v^*|^{2})}}{|v-v^*|}
\left|w(\ell ,\vartheta )d_2^i g_2 \right|^2dv^* dv\right)^{1/2},
\\
\le
Ce^{-\frac{q^\prime}{144}m_1^{2}}
\sqrt{
\left(\int e^{-\frac{q^\prime}{72}|v^*|^{2}}
\left|(w(\ell ,\vartheta )d_1^jg_1)_* \right|^2dv^*\right)
\left(\int e^{-\frac{q^\prime}{72}|v|^{2}}
\left|w(\ell ,\vartheta )d_2^i g_2 \right|^2dv\right)}.
\end{gather*}
By \eqref{hopdef} and Corollary \ref{lem:sigmaeigs}.2, we conclude
\begin{gather*}
\left|\langle {\bf 1}_{S_1} w^2K_{1}g_1, g_2\rangle\right|
\le
Ce^{-\frac{q^\prime}{144}m_1^{2}}
|g_1|_{\boldsymbol{\sigma},\vartheta} |g_2|_{\boldsymbol{\sigma},\vartheta}.
\end{gather*}
Since $m_1>0$ will be chosen large, this yields Theorem \ref{thm:upper} for $K_{1}$ restricted to $S_1$.  

Next fix $0<\tau<1$ and consider the region
$$
S_2=
\left\{\left|\frac{v-v^*}{|v-v^*|}\cdot v\right|\ge |v|^{\tau}\right\}
\cup
\left\{\left|\frac{v-v^*}{|v-v^*|}\cdot v^* \right|\ge |v^*|^{\tau}\right\}.
$$
Also $S_1^c=\{|v|< 2|v^*|<4|v|\}$.
Then on $S_1^c\cap S_2$, \eqref{bijbound00} yields
\begin{gather}
\left|w(\ell,\vartheta)(v)B_{ij}(v, v-v^*) \mu^{1/2}(v)\mu^{1/2}(v^*) \right|
\le
Cw(\ell,\vartheta)(v^*)\frac{e^{-q^{\prime\prime}(|v|^{2\tau}+|v^*|^{2\tau})}}{|v-v^*|},
\label{bijbound2}
\end{gather}
where $q^{\prime\prime}=\frac{q^\prime}{4}\frac{1}{2^{2\tau}}$.  Since $\tau>0$ we can use exactly the same estimates as in the previous case to establish
\begin{gather*}
\left|\langle {\bf 1}_{S_1^c\cap S_2} w^2(\ell,\vartheta)K_{1}g_1, g_2\rangle\right|
\le
Ce^{-q^{\prime\prime\prime}m_1^{2\tau}}
|g_1|_{\boldsymbol{\sigma},\vartheta} |g_2|_{\boldsymbol{\sigma},\vartheta},
\end{gather*}
with $q^{\prime\prime\prime}=q^{\prime\prime\prime}(\tau)>0$. And this grants Theorem \ref{thm:upper} for $K_{1}$ restricted to $S_1^c\cap S_2$ if $m_1>0$ is large enough.

It remains to estimate $K_{1}$ on $S_1^c\cap S_2^c$, where
$$
S_2^c=
\left\{\left|\frac{v-v^*}{|v-v^*|}\cdot v\right|\le |v|^{\tau}\right\}
\cap
\left\{\left|\frac{v-v^*}{|v-v^*|}\cdot v^* \right|\le |v^*|^{\tau}\right\},
$$
with  $0<\tau<1$. Over this region, we do not expect to get anymore general exponential decay  in 
$|v|$ and $|v^*|$ out of \eqref{bijbound0}-\eqref{bijbound00}.  However, fortunately,  Theorem \ref{lem:kernelestimate} allows us the possibility of finding polynomial decay in this region.

 Without loss of generality assume $|v|\ge 2^{1/2(1-\tau)}$, which means $\frac12 |v|^2-|v|^{2\tau}\ge 0$. Then, using \eqref{eq:vR}, \eqref{bijbound00} and $S_2^c$, with $0<\delta<1$ we have
\begin{gather*}
\left| w(\ell,\vartheta)(v)B_{ij}(v, v-v^*)\sqrt{\mu \mu_*}\right| 
\\
\le 
\frac{Cw(\ell,\vartheta)(v^*)}{|v-v^*|}
\left[1+{|v|^2-\left(\frac{(v-v^*)\cdot v}{|v-v^*|}\right)^2}\right]^{-(3+\delta)/2}
\\
\le 
\frac{Cw(\ell,\vartheta)(v^*)}{|v-v^*|}
\left[1+{|v|^2-|v|^{2\tau}}\right]^{-(3+\delta)/2}
\\
\le 
\frac{Cw(\ell,\vartheta)(v^*)}{|v-v^*|}
\left[1+\frac{1}{2}{|v|^2}\right]^{-(3+\delta)/2}.
\end{gather*}
In particular, since $|v|>2|v^*|$ on $S_1^c\cap S_2^c$ we have decay in both variables:
\begin{gather}
\left|w(\ell,\vartheta)(v) B(v, v-v^*)\sqrt{\mu_*\mu}\right| 
\le 
\frac{Cw(\ell,\vartheta)(v^*)}{|v-v^*|}
\frac{\left[1+{|v|^2}\right]^{-(3+\delta)/4}}{\left[1+{|v^*|^2}\right]^{(3+\delta)/4}}.
\label{bijbound3}
\end{gather}
And the next step is to use \eqref{bijbound3} to complete the estimate for $K_{1}$ on $S_1^c\cap S_2^c$.

By  \eqref{bijbound3} and the definition of $\varphi_1$ we have
\begin{gather*}
\left|\langle {\bf 1}_{S_1^c\cap S_2^c} w^2(\ell,\vartheta)K_1g_1, g_2\rangle\right|
\\
\le
C\int_{S_1^c\cap S_2^c} 
\frac{{\bf 1}_{\{|v|\ge m_1/3\}}}{|v-v^*|}
\frac{\left| (w(\ell,\vartheta)d_1^jg_1)_*\right|}{\left[1+{|v^*|}\right]^{(3+\delta)/2}}
\frac{\left|w(\ell,\vartheta)d_2^ig_2 \right|}{\left[1+{|v|}\right]^{(3+\delta)/2}}
 dv^* dv.
\end{gather*}
We next use Fubini and Cauchy-Schwartz to obtain
\begin{gather}
 =
C\int_{|v|\ge m_1/3}
\frac{\left|w(\ell,\vartheta)d_2^ig_2\right|}{\left[1+{|v|}\right]^{(3+\delta)/2}}
\left\{ \int_{v^*\in S_1^c\cap S_2^c} 
\frac{\left| (w(\ell,\vartheta)d_1^jg_1)_*\right| }{|v-v^*|\left[1+{|v^*|}\right]^{(3+\delta)/2}}dv^* \right\}dv
\nonumber
\\
\le
C
\left(\int_{|v|\ge m_1/3}
\left\{ \int_{S_1^c\cap S_2^c} 
\frac{\left| (w(\ell,\vartheta)d_1^jg_1)_*\right| }{|v-v^*|\left[1+{|v^*|}\right]^{(3+\delta)/2}} dv^* \right\}^2
 dv\right)^{1/2}
 \nonumber
 \\
 \times \left(\int_{|v|\ge m_1/3}w^2(\ell,\vartheta)(v)\frac{\left| d_2^ig_2\right|^2}{1+|v|^{3+\delta}} dv \right)^{1/2}.
 \label{???}
 \end{gather}
We now focus on bounding the term involving $(w(\ell,\vartheta)d_1^jg_1)_*$.
 
 The main difficulty with this term is controlling the singular factor $|v-v^*|^{-1}$.  The following is designed to control it.
First we remark that on $S_1^c\cap S_2^c$
$$
|v-v^*|=\left|\frac{(v-v^*)\cdot (v-v^*)}{|v-v^*|} \right|\le |v^*|^{\tau}+|v|^{\tau}\le (1+\sqrt{2})|v|^{\tau}.
$$ 
Therefore,
$
S_1^c\cap S_2^c \subset S_1^c\cap \left\{|v-v^*|\le 3 |v|^{\tau}\right\}.
$
Then
 \begin{gather*}
 \int_{v^*\in S_1^c\cap S_2^c} \frac{\left[1+{|v^*|}\right]^{-(3+\delta)/2}}{|v-v^*|}
\left| (w(\ell,\vartheta)d_1^jg_1)_*\right| dv^* 
\\
\le
\left( \int_{|v-v^*|\le 3 |v|^{\tau}} |v-v^*|^{-2}dv^*\right)^{1/2}
 \left( \int_{S_1^c\cap S_2^c} 
w^2(\ell,\vartheta)(v^*)\frac{\left| (d_1^jg_1)_*\right|^2}{\left[1+{|v^*|}\right]^{3+\delta}}
 dv^* \right)^{1/2}.
 \end{gather*}
  Moreover,
$$
\int_{ |v-v^*|\le 3|v|^{\tau}}|v-v^*|^{-2} dv^*\le C[1+ |v|]^{\tau}.
$$
 Combining these last two calculations yields
 \begin{gather*}
\int_{|v|\ge m_1/3}
\left\{ \int_{S_1^c\cap S_2^c} \frac{\left[1+{|v^*|}\right]^{-(3+\delta)/2}}{|v-v^*|}
\left| (w(\ell,\vartheta)d_1^jg_1)_*\right| dv^* \right\}^2
 dv
\\
\le 
 C
\left(\int_{S_1^c\cap S_2^c}{\bf 1}_{\{|v|\ge m_1/3\}}w^2(\ell,\vartheta)(v^*)
\frac{\left[1+{|v|}\right]^{\tau} }{1+{|v^*|}^{3+\delta}}
\left| (d_1^jg_1)_*\right|^2
 dv^*  dv\right)^{1/2}.
  \end{gather*}
  Since $v$ is comparable to $v^*$ on $S_1^c$, this is
  \begin{gather*}
\le 
 C
\left(\int_{S_1^c\cap S_2^c}{\bf 1}_{\{|v^*|\ge m_1/6\}}w^2(\ell,\vartheta)(v^*)
\frac{\left| (d_1^jg_1)_*\right|^2 }{1+{|v^*|}^{3+\delta-\tau}}
 dv^*  dv\right)^{1/2}.
  \end{gather*}
 And by Fubini's theorem, the integral in parenthesis is
   \begin{gather*}
=
\int_{|v^*|\ge m_1/6}w^2(\ell,\vartheta)(v^*)
\frac{\left| (d_1^jg_1)_*\right|^2 }
{1+{|v^*|}^{3+\delta-\tau}}
 \left( \int_{v\in S_1^c\cap S_2^c} dv \right)dv^*
 \\
 \le
C\int_{|v^*|\ge m_1/6}w^2(\ell,\vartheta)(v^*)
\frac{\left| (d_1^jg_1)_*\right|^2 }
{1+{|v^*|}^{3+\delta-\tau}}
 \left( \int_{|v-v^*|\le 3|v^*|^\tau} dv \right)dv^*
  \\
 \le
C
\int_{|v^*|\ge m_1/6}w^2(\ell,\vartheta)(v^*)
\frac{\left| (d_1^jg_1)_*\right|^2 }
{1+{|v^*|}^{3+\delta-\tau}}
\left[1+{|v^*|}\right]^{3\tau}dv^*.
  \end{gather*}
  Since $0<\delta, \tau<1$ are otherwise arbitrary, we can choose $\delta-4\tau$ to be arbitrarily close to $1$.  But this is not needed for the current case, so we merely choose 
  $\delta=\frac34$ 
  and 
  $\tau=\frac1{8}$ to obtain
   $\delta-4\tau=\frac14$.  
   Then, for $m_1>0$ large, the integral above is
\begin{gather*}
 \le
\frac{C}{m_1^{1/4}\log(m_1)}
\int_{\mathbb{R}^3}w^2(\ell,\vartheta)(v^*)
\frac{ \log(2+|v^*|)}{1+|v^*|^{3}}
\left| (d_1^jg_1)_*\right|^2dv^*.
\end{gather*}
We plug the result of the last few inequalities back into \eqref{???} to see that
\begin{gather*}
\left|\langle {\bf 1}_{S_1^c\cap S_2^c} w^2(\ell,\vartheta)K_1g_1, g_2\rangle\right|
 \le
\frac{C\left(\int_{\mathbb{R}^3}w^2(\ell,\vartheta)(v^*)\frac{ \log(2+|v^*|)}{1+|v^*|^{3}}\left| (d_1^jg_1)_*\right|^2dv^*\right)^{1/2}}{m_1^{1/8}\sqrt{\log(m_1)}}
\\
\times
\left(\int_{|v|\ge m_1/3}w^2(\ell,\vartheta)(v)\frac{\left| d_2^ig_2\right|^2}{1+|v|^{3+\delta}} dv \right)^{1/2}
\\
 \le
C
\frac{\left(\int_{\mathbb{R}^3}w^2(\ell,\vartheta)(v^*)\frac{ \log(2+|v^*|)}{1+|v^*|^{3}}\left| (d_1^jg_1)_*\right|^2dv^*\right)^{1/2}}{m_1^{1/2}\log(m_1)}
\\
\times
\left(\int_{\mathbb{R}^3}w^2(\ell,\vartheta)(v)\frac{\log(2+|v|)}{1+|v|^{3}}\left| d_2^ig_2\right|^2 dv \right)^{1/2}.
\end{gather*}
By \eqref{hopdef} and Corollary \ref{lem:sigmaeigs}.2 then
\begin{gather*}
\left|\langle {\bf 1}_{S_1^c\cap S_2^c} w^2(\ell,\vartheta)K_1g_1, g_2\rangle\right|
\le 
\frac{ C|g_1|_{\boldsymbol{\sigma},\vartheta}|g_2|_{\boldsymbol{\sigma},\vartheta}}{m_1^{1/2}\log(m_1)}.
\end{gather*}
This completes the estimate for $K_1$ with $m_1>0$ large.

Next define $K_2=K-K_1$.  In this definition of $K_2$ we have  a smooth cutoff function  
$0\le \varphi_2(r)=1-\varphi_1(r)\le 1$ 
 such that (for some $m_1>0$ large)
$\varphi_2(r)=1$ for $r\le m_1$ and $\varphi_2(r)=0$ if $r\ge 2m_1$.  
Then from \eqref{K12} and \eqref{hopdef}
\begin{gather*}
\langle w^2K_2g_1, g_2\rangle
=
\sum_{i,j}\int_{\mathbb{R}^3\times \mathbb{R}^3} \varphi_2(|v|+|v^*|)w^2B_{ij}(v,v-v^*)\sqrt{\mu_*\mu} 
(d_1^jg_1)_*
d_2^i g_2 dv^* dv,
\end{gather*}
where $w^2=w^2(\ell,\vartheta)$.  Theorem \ref{lem:kernelestimate} implies 
\[
\zeta_{ij}\equiv w^2(\ell,\vartheta)\varphi_2(|v|+|v^*|)  B_{ij}(v,v-v^*)\in 
L^2(\mathbb{R}^3\times \mathbb{R}^3). 
\]
Therefore, for any given $m_2>0,$ we can choose a $C_c^\infty $ function 
$
\varphi_2^{ij}(v,v^*)
$ 
such that 
\begin{gather*}
||\zeta_{ij}-\varphi_2^{ij}||_{L^2(\mathbb{R}^3\times \mathbb{R}^3)}\le \frac 1{m_2}, 
\\
{\rm supp}\{\varphi_2^{ij}\}\subset \{|v^*|+|v|\le C(m_2)<\infty\} .
\end{gather*}
We split $K_2$ into a ``small'' part and a ``compact'' part as
\begin{equation*}
\langle K_2g_1,g_2\rangle =
\langle K_{2c}g_1,g_2\rangle
+
\langle K_{2s}g_1,g_2\rangle.
\end{equation*}
Here we have used the splitting
$$
\zeta_{ij}=\varphi_2^{ij}+[\zeta_{ij}-\varphi_2^{ij}],
$$
to define
\begin{equation*}
\begin{split}
\langle K_{2c}g_1,g_2\rangle
&=
\sum_{i,j}
\int \varphi_2^{ij}(v,v^{*})
\sqrt{\mu_*\mu} (d_1^jg_1)_* d_2^i g_2 dv^* dv,
\\
\langle K_{2s}g_1,g_2\rangle
&=
\sum_{i,j}
\int [\zeta_{ij}-\varphi_2^{ij}]
\sqrt{\mu_*\mu} (d_1^jg_1)_* d_2^i g_2dv^* dv.
\end{split}
\end{equation*}
We introduce this second smooth cuttoff so that we can integrate by parts.
We estimate each of these terms separately.

The second term $\left| \langle K_{2s}g_1,g_2\rangle \right|$ is bounded by 
\begin{gather*}
\|\zeta_{ij}-\varphi_2^{ij}\|_{L^2(\mathbb{R}^3\times \mathbb{R}^3)} 
\|\sqrt{\mu_*\mu} (d_1^jg_1)_* d_2^i g_2||_{L^2(\mathbb{R}^3\times \mathbb{R}^3)} 
\\
\le 
\frac C{m_2}\left| \mu^{1/2}d_1^jg_1\right|_{0}
\left|\mu^{1/2}d_2^i g_2\right|_{0}
\le 
\frac C{m_2}\left| g_1\right| _{\boldsymbol{\sigma},\vartheta}
\left|g_2\right|_{\boldsymbol{\sigma},\vartheta},
\end{gather*}
where the last line follows from Corollary \ref{lem:sigmaeigs}.2.

After integrations by parts, the first term  is 
$$
\langle K_{2c}g_1,g_2\rangle =
\sum_{i,j}
\int \left(\underline{d}_1^j \underline{d}_2^i \{\varphi_2^{ij}(v,v^{*})\sqrt{\mu_*\mu}\} \right) 
g_1(v^{*})g_2(v) dv^* dv,
$$
where 
\begin{gather*}
\underline{d}_1^j
=
-\partial_{v^*_j}+\frac{v^*_j}{2}, 
~~~
\underline{d}_2^i
=
-\partial_{v_i}+\frac{v_i}{2}+w_1(v)v_i. 
\end{gather*}
Since $w_1(v)\le C$, we therefore have
$$
\left|\langle K_{2c}g_1,g_2\rangle \right|
\le 
C||\varphi_2^{ij}||_{C^2}\left\{ \int_{|v|\le C(m_2)}|
g_1|^2dv\right\} ^{1/2}\left\{ \int_{|v|\le C(m_2)}|g_2|^2dv\right\}^{1/2}.
$$
Combining the last two estimates we conclude that for any $m_2>0$
$$
|\langle w^2(\ell,\vartheta)K_2g_1,g_2\rangle |  
\le 
\frac C{m_2}
\left| g_1\right| _{\boldsymbol{\sigma},\vartheta}\left| g_2\right|_{\boldsymbol{\sigma},\vartheta}
+
C(m_2)|g_1 {\bf 1}_{C(m_2)}|_{0}  |g_2 {\bf 1}_{C(m_2)}|_{0}.
$$
We thus conclude Theorem \ref{thm:upper} for the $K_2$ part by first choosing $m_1>0$ large and then $m_2>0$ is chosen large enough.

Finally consider $|\langle w^2(\ell,\vartheta) \partial_i \sigma^i g_1,g_2\rangle |$.  For $m> 0$, we split the integration region into a compact part and a large part:
$$
\{|v|\le m\}\cup\{|v|\ge m\}.
$$
Then by \eqref{sigmai} and \eqref{destlambda}, 
\begin{eqnarray*}
\int_{|v|\ge m}  w^2(\ell ,\vartheta )\left|\partial_i \sigma^i g_1 g_2\right| dv
&\le& 
C \int_{|v|\ge m}w^2(\ell ,\vartheta ) \frac{\log\left(2+|v|\right)}{1+|v|^{3}} \left|g_1 g_2\right| dv,
\\
&\le&
\frac{C}{m^2} \int_{|v|\ge m} w^2(\ell ,\vartheta )\frac{\log\left(2+|v|\right)}{1+|v|} \left|g_1 g_2\right| dv.
\end{eqnarray*}
And then, using Corollary \ref{lem:sigmaeigs}.2 and Cauchy-Schwartz, the integral on the r.h.s. is
$
\le 
|g_1|_{\boldsymbol{\sigma},\vartheta}|g_2|_{\boldsymbol{\sigma},\vartheta}.
$
Thus,
\begin{equation*}
\int_{|v|\ge m}  w^2(\ell ,\vartheta )\left|\partial_i \sigma^i g_1 g_2\right| dv
\le 
\frac{C}{m^2} |g_1|_{\boldsymbol{\sigma},\vartheta}|g_2|_{\boldsymbol{\sigma},\vartheta}.
\end{equation*}
By Remark \ref{lem:sigmaeigs}.1, we see that $\sigma^i$ and $\partial_i \sigma^i$ have no finite singularities.  Therefore,
on $\{|v|\le m\}$, 
\eqref{sigmai} and \eqref{destlambda} imply that there is a constant $\tilde{C}(m)$ such that
$$
 \left|\partial_i \sigma^i \right|\le \tilde{C}(m).
$$
Therefore,
\begin{equation*}
\int_{|v|\le m} w^2(\ell ,\vartheta )\left|\partial_i \sigma^i g_1 g_2\right| dv
\le 
C(m) \left(\int_{|v|\le m}  \left|g_1\right|^2 dv\right)^{1/2}
\left(\int_{|v|\le m}  \left|g_1\right|^2 dv\right)^{1/2}.
\end{equation*}
We conclude that for any $m>0$
\begin{gather*}
|\langle w^2(\ell ,\vartheta )\partial_i \sigma^i g_1,g_2\rangle |  
\le 
\frac{C}{m^2}
\left| g_1\right| _{\boldsymbol{\sigma},\vartheta}\left| g_2\right|_{\boldsymbol{\sigma},\vartheta}
+
C(m)|g_1 {\bf 1}_{C(m)}|_{\vartheta}|g_2 {\bf 1}_{C(m)}|_{\vartheta}.
\end{gather*}
Since $m>0$ can be arbitrarily large, this completes the estimate for the last term $|\langle w^2(\ell ,\vartheta )\partial_i \sigma^i g_1,g_2\rangle|$  and thereby finishes the proof of Theorem \ref{thm:upper}.  
\qed \\

Next we will use Theorem \ref{thm:upper}  to deduce coercivity of the linear operator, $L$, via a standard compactness argument.

\medskip \noindent {\bf Corollary \ref{thm:upper}.1}.
{\it $\exists \delta >0$ such that 
$
\langle Lg,g\rangle \ge \delta |\{{\bf I-P}\}g|_{\boldsymbol{\sigma}} ^2.  
$
}

\begin{proof}  We use Theorem \ref{thm:upper} and the method of contradiction. In this case, we have a sequence of functions 
$\{g_n(v)\}_{n\ge 1}$ satisfying $|\{{\bf I-P}\} g_n|_{\boldsymbol{\sigma}}>0$ and
$$
\langle Lg_n, g_n \rangle < \frac{1}{n}  |\{{\bf I-P}\} g_n|_{\boldsymbol{\sigma}}^2.
$$
Without loss of generality suppose ${\bf P}g_{n}=0$ and $| g_n|_{\boldsymbol{\sigma}}\equiv 1$.  
Consider the inner product
$$
\langle g_1,g_2 \rangle_{\boldsymbol{\sigma}}=
\sum_{i,j}\int_{\mathbb{R}^3}\sigma^{ij}(v) \left( \partial_i g_1 \partial_j g_2 
+ \frac{v_i}{2}\frac{v_j}{2} g_1g_2\right) dv.
$$
Then there exists a $g_0(v)$ such that
$$
g_n\overset{\sigma}{\rightharpoonup}  g_0.
$$
In other words, $\langle g_n, h\rangle_{\boldsymbol{\sigma}}\to \langle g_0, h\rangle_{\boldsymbol{\sigma}}$ 
as $n\to\infty$ for all $h(v)$ which are bounded in $|\cdot|_{{\boldsymbol{\sigma}}}$.  By lower semi continuity,
$$
| g_0|_{\boldsymbol{\sigma}}\le 1.
$$
Equivalently
$
|g_0|_{\boldsymbol{\sigma}}^2 =\lim_{n\to\infty} \langle g_n, g_0\rangle 
\le \lim_{n\to\infty} |g_n|_{\boldsymbol{\sigma}} |g_0|_{\boldsymbol{\sigma}}=|g_0|_{\boldsymbol{\sigma}}.
$

From Lemma \ref{AKrepresent}  we can write 
\begin{equation}
\langle Lg_n,g_n\rangle =|g_n|_{\boldsymbol{\sigma}} ^2
-\langle \partial _i\sigma^ig_n,g_n\rangle -\langle Kg_n,g_n\rangle . 
\label{lEXPAND}
\end{equation}
We $claim$ that 
\begin{equation}
\begin{split}
\lim_{n\rightarrow \infty }\langle \partial _i\sigma ^ig_n,g_n\rangle
&=\langle \partial _i\sigma ^ig_0,g_0\rangle ,
\\
\lim_{n\rightarrow \infty
}\langle Kg_n,g_n\rangle &= \langle Kg_0,g_0\rangle . 
\end{split}
\label{ckLIMITS}
\end{equation}
These limits will follow directly from Theorem \ref{thm:upper}.  It is then a standard application of \eqref{ckLIMITS} to prove the coercivity.  

We first consider the limit of $\langle \partial _i\sigma ^ig_n,g_n\rangle$.  
Splitting
$$
\langle \partial _i\sigma ^ig_n,g_n\rangle
-\langle \partial _i\sigma ^ig_0,g_0\rangle 
=
\langle \partial _i\sigma ^i(g_n-g_0),g_n\rangle
+
\langle \partial _i\sigma ^ig_0,(g_n-g_0)\rangle. 
$$
From Theorem \ref{thm:upper} then
\begin{eqnarray*}
\left| \langle \partial _i\sigma ^ig_n,g_n\rangle
-\langle \partial _i\sigma ^ig_0,g_0\rangle \right|
&\le&
\frac {1}{8}\eta
\left| g_n-g_0\right| _{{\boldsymbol{\sigma}}}\left(\left| g_n\right|_{{\boldsymbol{\sigma}}}+\left| g_0\right|_{{\boldsymbol{\sigma}}}\right)
\\
&&+C(\eta)|(g_n-g_0) {\bf 1}_{C(\eta)}|_{0}\left(|g_n {\bf 1}_{C(\eta)}|_{0}+|g_0 {\bf 1}_{C(\eta)}|_{0}\right).
\end{eqnarray*}
Then by Corollary \ref{lem:sigmaeigs}.2,  
$|g_0 {\bf 1}_{C(\eta)}|_{0}\le C\left| g_0\right|_{{\boldsymbol{\sigma}}}$ 
and 
$\left| g_0\right|_{{\boldsymbol{\sigma}}}\le \left| g_n\right|_{{\boldsymbol{\sigma}}}=1$.
We thus have 
$$
\left| \langle \partial _i\sigma ^ig_n,g_n\rangle
-\langle \partial _i\sigma ^ig_0,g_0\rangle \right|
\le 
\frac {1}{2}\eta
+C(\eta)|(g_n-g_0) {\bf 1}_{C(\eta)}|_{0}.
$$
Similarly,
$$
\left| \langle Kg_n,g_n\rangle
-\langle Kg_0,g_0\rangle \right|
\le 
\frac {1}{2}\eta
+C(\eta)|(g_n-g_0) {\bf 1}_{C(\eta)}|_{0}.
$$
Furthermore $\partial _ig_n$ are bounded in $L^2\{|v|\le C(\eta)\}$ because $|g_n|_{\boldsymbol{\sigma}}
=1$; the Rellich-Kondrachov Compactness Theorem thus yields
$$
\lim_{n\uparrow\infty}|(g_n-g_0) {\bf 1}_{C(\eta)}|_{0}= 0. 
$$
Then \eqref{ckLIMITS} follows by first choosing $\eta$ small and second sending $n\uparrow \infty$.

By sending $n\uparrow \infty$ in \eqref{lEXPAND}, using \eqref{ckLIMITS}, we get
$$
0
=
1-\langle \partial_i \sigma^i g_0, g_0 \rangle
-\langle Kg_0, g_0 \rangle.
$$
We therefore conclude 
$$
\langle Lg_0, g_0 \rangle=
|g_0|_{\boldsymbol{\sigma}}^2-1\le 0.
$$
Since $L\ge 0$ we deduce $|g_0|_{\boldsymbol{\sigma}}=1$.     Equivalently, $\langle Lg_0, g_0 \rangle=0$.  We thus conclude that $g_0={\bf P} g_0$ by Lemma \ref{nullL}.  On the other hand, since ${\bf P}g_n=0$, we also conclude that $g_0=0$.  This contradicts $|g_0|_{\boldsymbol{\sigma}}=1$ and thereby establishes the result.
\end{proof}

Now that we have established Coercivity of the linearized Balescu-Lenard operator, we are ready to  prove Theorem \ref{thm:MAIN}.  Our arguments in this section are based on techniques developed in \cite{strainGUOed} combined with previous results in this work.  \\ 

\noindent {\it Proof of Theorem \ref{thm:MAIN}.}
We will prove that a solution to the linearized Balescu-Lenard equation satisfying the assumptions of Theorem \ref{thm:MAIN} also satisfies:
\begin{equation}
\label{energyBOUND}
\frac{d}{dt} |f|^2_{\vartheta}(t)+| f|^2_{\boldsymbol{\sigma},\vartheta}(t)\le 0,
\end{equation}
Here we use the notation $|f|^2_{\vartheta}$ and $| f|^2_{\boldsymbol{\sigma},\vartheta}$ loosely in the sense that we establish the above only up to equivalent norms.
Then we will show that \eqref{energyBOUND} implies 
\eqref{timeDECAY}.  

We prove \eqref{energyBOUND}.  Multiply \eqref{linearizedBL} by $f$ and integrate over $v\in\mathbb{R}^3$ to obtain
$$
\frac{d}{dt} |f|^2_0+\langle Lf, f\rangle =0.  
$$
We establish \eqref{energyBOUND} for $\ell=\vartheta=0$ by Corollary \ref{thm:upper}.1.  

Next assume $\ell\in\mathbb{R}$ and $0\le \vartheta\le 2$.  In this case we multiply \eqref{linearizedBL} by $w^2(\ell,\theta) f$
and integrate over $v\in\mathbb{R}^3$ to obtain
$$
\frac{d}{dt} |f|^2_{\vartheta}+\langle w^2(\ell,\theta) Lf, f\rangle =0.  
$$
Our goal is to get a lower bound for $\langle w^2(\ell,\theta) Lf, f\rangle$.  By Lemma \ref{AKrepresent} we can write
\begin{equation}
\langle w^2(\ell,\vartheta) Lg, g\rangle=
\left|g\right|_{\boldsymbol{\sigma},\vartheta}^2 
+ \langle \partial_i (w^{2})\sigma^{ij} \partial_j g, g\rangle 
- \langle w^{2}\partial_i\sigma^i g, g\rangle 
- \langle w^{2}K g, g\rangle.
\label{claimLINEARexp}
\end{equation}
Recall that ${\bf 1}_{C}$ is the indicator function of the set $\{|v|\le C\}$.
We {\it claim} that 
\begin{equation}
\langle w^{2}(\ell,\vartheta) [L g], g\rangle \ge \delta_q
\left| g\right|_{\boldsymbol{\sigma},\vartheta}^2 
-
C(\eta) \left|{\bf 1}_{C(\eta)} g\right|_{0}^2,
\label{claimLINEAR}
\end{equation}
where $\delta_q>0$ depends on $q$ from \eqref{weight} and $\eta>0$ chosen small enough.  Notice that by  plugging \eqref{claimLINEAR} into the differential equality at the beginning of this paragraph and adding the result to \eqref{energyBOUND} for the case $\vartheta=\ell=0$ yields \eqref{energyBOUND} in the general case (at least up to an equivalent norm).  

To establish \eqref{energyBOUND}, it remains to prove  \eqref{claimLINEAR}.  We do this first assuming
$\ell\in\mathbb{R}$ and $0\le \vartheta< 2$.  Later we will handle the case
$\vartheta = 2$ separately.
We estimate each of the last three terms on the right side of \eqref{claimLINEARexp}.  By Lemma \ref{lem:sigma} and 
\eqref{wDERIVATIVE} we can write 
\begin{equation*}
\partial_i (w^{2}(\ell,\vartheta)(v))\sigma^{ij}(v)=w^2(v) w_1(v) \lambda_1(v)v_j.
\end{equation*}
From the definition of $w_1(v)$ in \eqref{wDERIVATIVE}, 
$\left| w_1(v)\right| \le C(1+|v|^2)^{\frac{\vartheta}{2}-1}$.
By Lemma \ref{lem:sigmaeigs}, 
$\left|\lambda_1(v)v_j \right| \le C \frac{\log(2+|v|)}{1+|v|^2}$. Thus for any $m^\prime>0$ 
\begin{gather*}
\begin{split}
\left| \langle \partial_i (w^{2})\sigma^{ij} \partial_j g, g\rangle \right|
&\le C 
\int \frac{w^2(\ell,\vartheta)\log(2+|v|)}{1+|v|^{2}}
\frac{\left|\partial_j g \right| \left| g \right|}{[1+|v|]^{\frac{\vartheta}{2}-1}}dv.
\end{split}
\end{gather*}
For $m>0$ large,  we split the integral into a bounded part and an unbounded part.  We have
\begin{gather*}
\begin{split}
\int_{|v|> m} w^2\frac{\log(2+|v|)}{1+|v|^{2}}
\frac{\left|\partial_j g \right| \left| g \right|}{[1+|v|]^{\frac{\vartheta}{2}-1}}dv
&\le 
C m^{\frac{\vartheta-2}{2}}\int_{|v|> m}w^2 \frac{\log(2+|v|)}{1+|v|^{2}}
\left|\partial_j g \right| \left| g \right|dv
\\
&\le 
C m^{\frac{\vartheta-2}{2}} | g|^2_{\boldsymbol{\sigma},\vartheta}
\le
\frac{\eta}{3} | g|^2_{\boldsymbol{\sigma},\vartheta},
\end{split}
\end{gather*}
where we have used Cauchy-Schwartz, Corollary \ref{lem:sigmaeigs}.2, $0\le \vartheta<2$ and $m>0$ chosen large enough.  
We use Cauchy-Schwartz on the bounded part to obtain
\begin{gather*}
\begin{split}
\int_{|v|\le m} w^2\frac{\log(2+|v|)}{1+|v|^{2}}
\frac{\left|\partial_j g \right| \left| g \right|}{[1+|v|]^{\frac{\vartheta}{2}-1}}dv
&\le 
C(m) | g|_{\boldsymbol{\sigma}}\left|{\bf 1}_{m} g\right|_{0}, \\
&\le 
\frac{\eta}{3} | g|_{\boldsymbol{\sigma},\vartheta}^2+C(\eta)\left|{\bf 1}_{m} g\right|_{0}^2.
\end{split}
\end{gather*}
We estimate the last two terms on the right side of \eqref{claimLINEARexp} in the same way using
Theorem \ref{thm:upper}.
This establishes the {\it claim} for $0\le \vartheta<2$.

For $\vartheta=2$, to prove \eqref{claimLINEARexp} we will split linear operator in a different way.    
Specifically, split $L=-A-K$ (as in Lemma \ref{AKrepresent}) and define 
$M(v)\equiv \exp\left(\frac{q}{4}\left(1+|v|^2\right)\right)$. First we can show that there is $\delta_q>0$
such that 
\begin{equation*}
\begin{split}
-\langle w^{2}(\ell,2) [A g], g\rangle 
\ge &\delta_q 
\int_{\mathbb{R}^3} (1+|v|^{2})^{ \ell }\sigma^{ij}(v) \left( \partial_i [Mg] \partial_j [Mg] + \frac{v_i}{2}\frac{v_j}{2}|Mg|^2\right) dv
\\
&-
C(\delta_q) \left|{\bf 1}_{C(\delta_q)} g\right|_{0}^2.  
\end{split}
\end{equation*}
Second we can establish 
\begin{equation*}
\int_{\mathbb{R}^3} (1+|v|^{2})^{ \ell }\sigma^{ij}(v) \left( \partial_i [Mg] \partial_j [Mg] + \frac{v_i}{2}\frac{v_j}{2}|Mg|^2\right) dv
\ge 
\delta_q | g|_{\boldsymbol{\sigma},\vartheta}^2 - C(\delta_q)\left|{\bf 1}_{C(\delta_q)} g\right|_{0}^2, 
\end{equation*}
where $\delta_q=1-q^2-\frac{\eta}{2}>0$ since $0<q<1$ and $\eta>0$ can be chosen
arbitrarily small. This is enough to establish \eqref{claimLINEAR}
because the $K$ part is controlled by Theorem \ref{thm:upper}.
Using slightly different notation, this exact result was shown for the linearized Landau collision operator in \cite[Lemma 9]{strainGUOed} in equations (64) and (65) of that paper.  Since the proof is very much the same in this case, we will not repeat it.

Now that we have established our {\it claim} and thereby the differential inequality \eqref{energyBOUND}, it remains to show that this implies exponential decay.  As in other works on time decay problems for soft potentials \cite{MR575897,strainGUOed}, a key point is to split $|f|_\vartheta^2(t)$ into a time dependent low velocity part 
\begin{equation*}
E = \{1+|v|\le \rho (1+t)^{p^\prime}\},
\end{equation*}
Then from  Corollary \ref{lem:sigmaeigs}.2, we have 
\begin{equation*}
|f|^2_{\boldsymbol{\sigma}}(t) \ge \frac{C}{\rho(1+ t)^{p^\prime}} |f{\bf 1}_E|^2_{0}(t),  
\end{equation*}
where ${\bf 1}_E$ is just the indicator of the set $E$.
Notice that we have ignored the weak logarithmic factor.  
Plugging this into the the differential inequality \eqref{energyBOUND} we
obtain 
\begin{equation*}
\frac{d}{dt}|f|_{0}^2(t)+\frac{C}{\rho(1+ t)^{p^\prime}} |f {\bf 1}_E|^2_{0}(t)\le 0.
\end{equation*}
Thus
\begin{equation*}
\frac{d}{dt}|f|_{0}^2(t)+\frac{C}{\rho(1+ t)^{p^\prime}} |f|_{0}^2(t)
\le 
\frac{C}{\rho(1+ t)^{p^\prime}} |f {\bf 1}_{E^c}|^2_{0}(t).
\end{equation*}
Define $\lambda=\frac{C}{\rho p}$ where for now $p=1- p^\prime$ and 
$p^\prime>0$ is otherwise arbitrary. Then 
\begin{equation*}
\frac{d}{dt}\left(e^{\lambda (1+t)^{p}}|f|_{0}^2(t)\right) 
\le
\lambda p (1+ t)^{p-1}  e^{\lambda (1+t)^{p}}|f {\bf 1}_{E^c}|^2_{0}(t).
\end{equation*}
The integrated form is 
\begin{equation*}
|f|_{0}^2(t) \leq e^{-\lambda (1+t)^{p}}|f_0|_{0}^2(t)+ \lambda
p e^{-\lambda (1+t)^{p}} 
\int_0^t (1+s)^{p-1}e^{\lambda (1+s)^{p}}|f {\bf 1}_{E^c}|^2_{0}(t) ds.
\end{equation*}
Since $|f {\bf 1}_{E^c}|^2_{0}(t)$
is on 
$
E^c = \{1+|v|> \rho (1+t)^{p^\prime}\}
$ 
we have
\begin{equation*}
|f {\bf 1}_{E^c}|^2_{0}(s)  \le C e^{-\frac{q}{2}\rho^\vartheta (1+s)^{\vartheta p^\prime}}|f|^2_{\vartheta}(s) .
\end{equation*}
In the last display we have used the region and 
\begin{equation*}
1\le \exp\left(\frac{q}{2}(1+|v|^2)^{\frac{\vartheta}{2}}-\frac{q}{2}(1+|v|)^{\vartheta}\right)
 \le 
 \exp\left(\frac{q}{2}(1+|v|^2)^{\frac{\vartheta}{2}%
}\right) e^{-\frac{q}{2}\rho^\vartheta (1+s)^{\vartheta p^\prime}}.
\end{equation*}
The integrated form of \eqref{energyBOUND} implies 
\begin{gather}
|f|_{0}^2(t) \leq e^{-\lambda t^{p}} \left( |f_0|_{0}^2(t)
+ \lambda p |f_0|_{\vartheta}^2 \int_0^t (1+s)^{p-1} 
e^{\lambda (1+s)^{p}- \frac{q}{2}\rho^\vartheta (1+s)^{\vartheta p^\prime}} ds\right).
\label{finite?}
\end{gather}
The biggest exponent $p$ that we can allow with this splitting is $%
p=\vartheta p^\prime$; since also $p=-p^\prime+1$ we have $p^\prime=%
\frac{1}{\vartheta+1} $ so that 
\begin{equation*}
p=-\frac{1}{\vartheta+1}+1=\frac{\vartheta}{\vartheta+1}.
\end{equation*}
Further choose $\rho>0$ large enough so that 
$\lambda=\frac{C}{\rho p}<\frac{q}{2}\rho^\vartheta$.  Therefore 
the right side of \eqref{finite?} is finite.
This completes the proof of decay. \qed \\

\noindent {\bf Acknowledgements}.  The author would like to express his gratitude to Yan Guo for suggesting that he study this equation.  He also thanks Cl{\'e}ment Mouhot for several stimulating discussions regarding this work.  This work was supported by an NSF Mathematical Sciences Postdoctoral Research Fellowship.

\end{document}